\documentclass[a4paper,oneside,10pt]{article}
\usepackage{amsmath}
\usepackage{amssymb}
\usepackage{amsfonts}
\usepackage{graphicx}
\usepackage{color}
\usepackage{amsthm}
\usepackage{cite}
\usepackage{enumitem}
\usepackage{hyperref} 

\newtheoremstyle{uprightstyle}
{3pt} 
{3pt} 
{\upshape} 
{} 
{\bfseries} 
{.} 
{ }
{}

\theoremstyle{uprightstyle}

\newtheorem{theorem}{Theorem}[section]

\newtheorem{lemma}[theorem]{Lemma}

\newtheorem{proposition}[theorem]{Proposition}
\newtheorem{remark}[theorem]{Remark}

\newenvironment{proof}[1][Proof]{\noindent \textbf{#1.} }{\ \ $\Box$}

\numberwithin{equation}{section}

\renewcommand \arraystretch{1.5}
\begin{document}
	\title{BSDEs driven by \texorpdfstring{$ G $}{}-Brownian motion with time-varying uniformly continuous generators}
	\author{Bingru Zhao \thanks{Zhongtai Securities Institute for Financial Studies,
			Shandong University, Jinan, Shandong 250100, China. bingruzhao@mail.sdu.edu.cn. Research supported by NSF (No. 12326603, 11671231). }}
	\date{} 
	\maketitle
	
	\textbf{Abstract}. In this paper, we study the backward stochastic differential
	equations driven by $G $-Brownian motion under the condition that the
	generator is time-varying Lipschitz continuous with respect to $y $ and time-varying  uniformly continuous with respect to $z $. With the help of linearization method and the $ G $-stochastic analysis techniques, we construct the approximating sequences of $ G $-BSDE and obtain some precise a priori estimates. By combining this with the approximation method, we prove the existence and uniqueness of the solution under the time-varying conditions, as well as the comparison theorem.	
	
	{\textbf{Key words}. } $ G $-expection, $ G $-Brownian motion, BSDE, time-varying,
	uniformly continous
	
	\textbf{AMS subject classifications.} 93E20, 60H10, 35K15
	
	\addcontentsline{toc}{section}{\hspace*{1.8em}Abstract}
	\section{Introduction}

	In 1990, Pardoux and Peng \cite{pardoux1990adapted} first proposed backward stochastic differential equations (BSDE) in a Wiener probability space $(\Omega, \mathcal{F},P) $ on finite interval $ [0,T] $:
	\begin{equation}
		Y_{t}=\xi +\int_{t}^{T} g(s,Y_{s},Z_{s})ds-\int_{t}^{T}Z_{s}dB_{s}
	\end{equation}
	and originally constructed the fundamental framework of BSDE under nonlinear cases, where $\Omega$ is the space of all continuous paths and $P$ is the Wiener measure. And it demonstrates the existence and uniqueness of solution for BSDE with Lipschitz continuous generators, which leads to the emergence of research on the well-posedness theory of nonlinear BSDE. This study laid the foundation for BSDE theory, establishing them as a crucial research tool in the fields of partial differential equations(PDE), stochastic control and mathematical finance\cite{buckdahn2010probabilistic,coquet2002filtration,el1997backward}.
	
	Motivated by financial problems with model uncertainty and based on a probabilistic interpretation of fully nonlinear PDE, Peng\cite{peng2005nonlinear,peng2007g,peng2010nonlinear} non-trivially extended the classical BSDE framework and constructed a time-consistent sublinear expectation framework, called $G $-expectation ($\hat{\mathbb{E}}[\cdot] $), which is more suitable for fully nonlinear cases without the need for the probability framework. The corresponding regular
	process $(B_t)_{t \geq 0} $ is called G-Brownian motion and the stochastic calculus of $ \text {Itô's} $ type is also established in Peng\cite{peng2008multi}. Thus, due to the $ G $-martingale representation theorem, the backward stochastic differential equations driven by $G $-Brownian motion ($ G $-BSDE) has the following form:
	\begin{equation}  \label{eq17}
		Y_{t}=\xi+\int_{t}^{T} f\left(s, Y_{s}, Z_{s}\right) d s+\sum_{i, j=1}^{d}
		\int_{t}^{T} g_{i j}\left(s, Y_{s}, Z_{s}\right) d\left\langle B^{i},
		B^{j}\right\rangle_{s}-\int_{t}^{T} Z_{s} d B_{s}-\left(K_{T}-K_{t}\right),
	\end{equation}
	where $K $ is a non-increasing and continuous $G $-martingale.
	
	By comparing the classical BSDE with $G $-BSDE, we can observe that due to the existence of non-increasing $G $-martingale, the methods used in the classical BSDE framework cannot be applied in the $G$-framework\cite{chen2000infinite}. Consequently, many studies have been devoted to overcoming this difficulty, which makes the study of $ G $-BSDE a challenging and interesting problem.
	
	For instance, with the help of fully nonlinear parabolic PDE and Galerkin-type approximation, Hu et al. \cite{hu2014backward} obtained the existence and uniqueness of $ G $-BSDE with Lipschitz generators. Liu\cite{liu2020multi} further explored multi-dimensional $ G $-BSDEs and solved the well-posedness problem under these conditions. Wang and Zheng \cite{wang2021backward} considered the G-BSDE under the condition that the generators \( f \) and \( g \) are uniformly continuous with respect to the variable \( z \). Following this, inspired by the research results related to BSDE, Hu et al. \cite{hu2020bsdes} relaxed the assumptions and considered the $G $-BSDE with time-varying Lipschitz generators, where the time-varying coefficients need to satisfy certain integrability condition:
	\begin{equation}\label{eq25}
		\int_{0}^{T}u(s)+v^2(s)ds < \infty. 
	\end{equation}
	In addition, by using a new type of nonlinear expectation, Soner et al. \cite{soner2012wellposedness,soner2013dual} proposed a second-order backward stochastic differential equation(2-BSDE). Lin et al. \cite{lin2bsde2020} established the well-posedness of 2-BSDE with random terminal time. Moreover, the research on quadratic G-BSDEs has presented significant results, and readers can refer to references \cite{hu2020quadraticbsde,huQuadraticGBSDEs2022} for further information.
	
	In this paper, we would like to study the $ G $-BSDE with time-varying uniformly continuous and Lipschitz continuous generators, which involves the integrability condition \eqref{eq25}. The main idea is to construct the supremum approximation equation and the infimum approximation equation for $ G $-BSDE \eqref{eq17} (see \eqref{eq5}, \eqref{eq6}), based on the fact that time-varying uniformly continuous generators and Lipschitz generators can be approximated by a sequence of Lipschitz generators. Combining this with linearization techniques and nonlinear stochastic analysis knowledge, this paper aims to deduce the existence and uniqueness of solution and comparison theorem. This idea differs from the principle of contraction mapping in classical BSDE, as its core lies in constructing the approximating sequences of the generator. Therefore, the convergence of approximating sequences of $ G $-BSDE is an crucial issue that needs to be addressed. 
	
	Since the time-varying coefficients are not bounded functions, the main difficulty is ensuring that the time-varying coefficients do not affect the convergence of approximating sequences of the $ G $-BSDE \eqref{eq17}. Additionally, selecting a linear decomposition function which is suitable for the time-varying conditions is also challenging. To overcome these difficulties, we fully utilize nonlinear stochastic analysis techniques and the properties of the approximation sequences of the generator, allowing us to obtain more precise estimates for $Y $, $Z $ and the difference between the approximating sequences of the $ G $-BSDE. This approach effectively prevents the blow up of time-varying coefficients and enables us to obtain the main results of the paper. 
	
	The paper is organized as follows: In section 2, we briefly introduce the relevant preliminary knowledge and important theorems under the $G $-expection framework; In Section 3, we study the impotant results for $ G $-BSDEs under time-varying uniformly continuous and Lipschitz continuous.
	
	\section{Preliminaries}
	
	This section will introduce some basic notions of $G $-expectation theory
	and important results of $G $-BSDEs with time-varying Lipschitz condition.
	These will be used in the following sections. More details can be found in 
	\cite{peng2007g,peng2008multi,peng2010nonlinear}.
	\subsection{\texorpdfstring{$ G $}{}-expection}
	Denote $\Omega = C_0^d(\mathbb{R}^{+})$ as the space consisting of all $\mathbb{R}^{d}$-valued continuous path functions $(\omega_t)_{t\geq 0}$
	originating at 0, where the metric distance is  
	\begin{equation*}
		\rho\left(\omega^{1}, \omega^{2}\right):=\sum_{n=1}^{\infty} \frac{1}{2^{n}}%
		\left[\left(\max _{t \in[0, n]}\left|\omega_{t}^{1}-\omega_{t}^{2}\right|%
		\right) \wedge 1\right],\quad \omega_{t}^{1},\omega_{t}^{2}\in \Omega .
	\end{equation*}
	Then, it is clear that $(\Omega,\rho)$ is a complete metric space. Set $\Omega_T := \left\{\omega_{\cdot \wedge T }:\omega \in \Omega \right\}$ for $T\in [0,\infty)$ and $\mathcal{B}(\Omega)$ is the Borel $\sigma$-algebra of $\Omega$. The $d $-dimensional canonical process $(B_t)_{t \geq 0}$ is defined as follows. For
	any $t\in [0,T]$, $\omega \in \Omega$, $B_t(\omega) := \omega_t$, which is
	called the $G $-Brownian motion. In particular, $B_1$ is $G $-normally
	distributed and $B_t \overset{d}{= }\sqrt{t}B_1. $ Let $t_1$
	$ ,..., $ $t_k$ is a partition on $[0,T]$, define 
	\begin{equation*}
		L_{ip}\left(\Omega_{t}\right):=\left\{\varphi\left(B_{t_{1}},
		\ldots, B_{t_{k}}\right): k \in \mathbb{N}, t_{1}<\cdots<t_{k} \in[%
		0, t], \varphi \in C_{b.L_{ip} }\left(\mathbb{R}^{k \times
			d}\right)\right\}
	\end{equation*}
	and 
	\[
	L_{ip}(\Omega):=\bigcup_{t \geq 0} L_{ip}\left(\Omega_{t}\right),
	\]
	where $C_{b.L_{ip} }\left(\mathbb{R}^{k \times d}\right)$ is the space
	consisting of all $\mathbb{R}^{k \times d}$-valued bounded and Lipschitz continuous functions.
	
	By the nonlinear parabolic PDE, Peng \cite{peng2010nonlinear} constructed the $G $-expectation space $(\Omega,L_{ip}(\Omega),\hat{\mathbb{E}}[\cdot],\hat{\mathbb{E}}_t[\cdot]_{t\geq 0})$, where the monotone and sublinear function $G:\mathbb{S}_d \rightarrow \mathbb{R} $ satisfies that 
	\begin{equation*}
		G(A)=\frac{1}{2}\hat{\mathbb{E}}[\langle AB_1,B_1 \rangle], \quad
		\forall A \in \mathbb{S}_d.  
	\end{equation*}
	We define $\mathbb{S}_d $ as the space of all $d\times d $ symmetric matrices, and $\mathbb{S}_d^{+} $ as the space of all $d\times d $ symmetric positive definite matrices. This paper assumes $G$ is non-degenerate, implying that there exists two constants $0< \underline{\sigma}^2 \leq \bar{\sigma}^2 <
	\infty$ such that  
	\begin{equation*}
		\frac{1}{2}\underline{\sigma}^2\text{tr}[A_1-A_2] \leq G(A_1)-G(A_2) \leq \frac{1}{2}\bar{\sigma}^2\text{tr}[A_1-A_2],\quad \forall A_1 \geq A_2.  
	\end{equation*}
	Thus, there is a subset $\Gamma \subset \mathbb{S}_d^{+}$ that is bounded, convex and closed, such that
	\begin{equation}\label{eq26}
		G(A) = \frac{1}{2}\underset{B\in \Gamma}\sup \text{tr}[AB].
	\end{equation}
	The following is the representation theorem for 
	a sublinear expectation, along with important results of $ G $-BSDEs, which will be widely used throughout the paper. To begin with, we first show that the $G $-expectation space required for this paper. For any $p \geq 1$ and fixed $T \geq 0$, 
	\begin{itemize}[noitemsep, leftmargin=*]
		\item $L_{G}^{p}\left(\Omega_{T}\right)$:= the completion of $L_{i
			p}\left(\Omega_{T}\right)$ under the norm $\|X\|_{L^p_G(\Omega_{T})} = \hat{%
			\mathbb{E}}\left[|X|^p\right ]^{1 / p}$(resp. $L_{G}^{p}\left(\Omega
		\right)$);  
		\item $M_{G}^{0}(0, T):=\left\{\eta_{t}=\sum_{j=0}^{N-1} \xi_{j} I_{\left[%
			t_{j}, t_{j+1}\right)}(t): 0=t_{0}<\cdots<t_{N}=T, \xi_{i} \in
		L_{ip}\left(\Omega_{t_{i}}\right)\right\}$; 
		\item $M_{G}^{p}(0, T):= $ the completion of $M_{G}^{0}(0, T)$ under the
		norm $\|X\|_{M_{G}^{p}} = \left\{\hat{\mathbb{E}}[\int_{0}^{T}\left|X_{s}%
		\right|^{p} d s]\right\}^{1 / p} $; 
		
		\item $H_{G}^{p}(0, T)$:= the completion of $M_{G}^{0}(0, T) $ under the
		norm $\|X\|_{H_{G}^{p}}=\left\{\hat{\mathbb{E}}[(\int_{0}^{T}\left|X_{s}%
		\right|^{2} d s)^{p / 2}]\right\}^{1 / p} $; 
		
		\item $S_{G}^{0}(0, T):=\left\{h \left(t, B_{t_{1} \wedge t}, \ldots,
		B_{t_{k} \wedge t}\right): k \in \mathbb{N}, t_{1}<\cdots<t_{k} \in[0, T], h
		\in C_{b \cdot {\ L_{ip} }}\left(\mathbb{R}^{n+1}\right)\right\}$; 
		
		\item $S_{G}^{p}(0, T):= $ the completion of $S_{G}^{0}(0, T)$ under the
		norm $\|X\|_{S_{G}^{p}} = \left\{\hat{\mathbb{E}}\left[\sup _{t \in
			[0,T]}|X_t|^p\right]\right\}^{1 / p} $; 
		
		\item $\mathfrak{S}_{G}^{P}(0, T):$= the triplet of processes $(Y, Z, K) $
		satisfying that $Y \in S_{G}^{p}(0, T)$ , $Z \in H_{G}^{p}(0, T)$, $K_{T} \in
		L_{G}^{p}\left(\Omega_{T}\right)$ and $K $ is a
		decreasing $G $-martingale with $K_{0}=0 $. 
	\end{itemize}
	
	\begin{theorem}
		\textnormal{\cite{denis2011function,hu2014backward}}  Let $\mathcal{P}$ is the set of	probability measures on $(\Omega,\mathcal{B}(\Omega))$. Then, there exists a tight set $\mathcal{P}$ such that for all $	\xi \in L_G^1(\Omega)$,  
		\begin{equation*}
			\hat{\mathbb{E}}[\xi] = \underset{P \in \mathcal{P}}\sup E_P[\xi] = \underset%
			{P \in \mathcal{P}}\max E_P[\xi].  
		\end{equation*}
	\end{theorem}
	\noindent Therefore, by the tight set $\mathcal{P}$, we define capacity as 
	\begin{equation*}
		c(A) := \sup_{P \in \mathcal{P}} P(A), \quad \forall A \in \mathcal{B}(\Omega). 
	\end{equation*}
	A set $A \in \mathcal{B}(\Omega) $ is called polar if it satisfies that $c(A) = 0$. A property is said to hold "quasi-surely"(q.s.) if it holds outside a polar set. In the following, two random variables $X$ and $Y$ are considered equivalent if $X = Y$ q.s. 	
	\begin{proposition}
		\textnormal{\cite{hu2014comparison}} Let $\alpha \geq 1$, $p \in(0, \alpha] $, if $\eta \in H_{G}^{p}(0, T) $,
		then\newline
		\begin{equation*}
			\underline{\sigma}^{p} c_{p} \hat{\mathbf{E}}_{t}\left[\left(\int_{t}^{T}%
			\left|\eta_{s}\right|^{2} d s\right)^{p / 2}\right]  \leq \hat{\mathbf{E}}%
			_{t}\left[\sup _{u \in[t, T]}\left|\int_{t}^{u} \eta_{s} d B_{s}\right|^{p}%
			\right] \leq \bar{\sigma}^{p} c_{p} \hat{\mathbf{E}}_{t}\left[%
			\left(\int_{t}^{T}\left|\eta_{s}\right|^{2} d s\right)^{p / 2}\right] , 
		\end{equation*}
		where $\underline{\sigma}^{2}:=-\hat{\mathbf{E}}\left(-B_{1}^{2}\right)<\hat{%
			\mathbf{E}}\left(B_{1}^{2}\right):=\bar{\sigma}^{2} $ and $c_{p}$ is a
		constant depending on $p$. 
	\end{proposition}
	\subsection{\texorpdfstring{$ G $}{}-BSDE with time-varying Lipschitz generators}
	Next, we will introduce the important results of $G $-BSDE with time-varying generators, as well as linear $G $-BSDE. Hu et al. \cite{hu2020bsdes} proved  the
	existence and uniqueness of solution for $G $-BSDE \eqref{eq17} with the time-varying Lipschitz condition, where the generators  
	\begin{equation*}
		f(t, \omega, y, z),\text{ }g(t, \omega, y, z):[0, T] \times \Omega \times 
		\mathbb{R} \times \mathbb{R}^{d} \rightarrow \mathbb{R}  
	\end{equation*}
	and the terminal $\xi $ satisfy (A1)-(A3).\newline
	(A1) There exists a constant $\beta \geq 0 $ such that
	for any $\left (y, z\right) \in \mathbb{R} \times \mathbb{R}^{d} $, $%
	f(\cdot, \cdot, y, z)$, $g(\cdot, \cdot, y, z) \in M_{G}^{1}(0, T) $ and  
	\begin{equation*}
		\int_{0}^{T}(|f(s, 0,0)|+\sum_{i, j=1}^{d}\left|g_{i j}(s, 0,0)\right|)ds
		\in L_{G}^{2+\beta}\left(\Omega_{T}\right) .  
	\end{equation*}
	(A2) There exist processes $%
	u(t) $ and $v(t) $ that is positive and non-random function, such that  
	\begin{equation*}
		\left|\ell(t, y, z)-\ell\left(t, y^{\prime}, z^{\prime}\right)\right| \leq
		u(t)\left|y-y^{\prime}\right|+v(t)\left|z-z^{\prime}\right|,  
	\end{equation*}
	where $\ell=f$, $ g_{i j} $ .  \newline
	(A3) The $ u(t) $ and $ v(t) $ satisfy the condition of square integrability:  
	\begin{equation*}
		\Lambda^{T}(u, v):=\int_{0}^{T}\left(u(s)+v^{2}(s)\right)dt<\infty.  
	\end{equation*}
	
	\begin{theorem}
		Suppose that $f $, $g $
		satisfy the assumptions $ \textnormal{(A1)-(A3)} $ and $\xi \in L_{G}^{2+\beta}\left(\Omega_{T}\right) $ for some $\beta>0 $. Then the $G $-BSDE \eqref{eq17} admits a unique solution 
		$(Y, Z, K) \in \mathfrak{S}_{G}^{2}(0, T) $. 
	\end{theorem}
	
	\begin{proposition}
		Given that $f^l $, $g^l $
		satisfy the assumptions $ \textnormal{(A1)-(A3)} $ and $\xi^l \in L_{G}^{2+\beta}\left(\Omega_{T}\right) $ for some $\beta>0 $. And suppose that $\left(Y^{l}, Z^{l}, K^{l}\right) \in \mathfrak{S}_{G}^{2}(0, T)$, $l=1,2 $ is the solution of $G $-BSDE \eqref{eq17} associated with	the data $\left(\xi^{l} , f^{l}, g^{l}\right)$. Then there exists a constant 
		$C(\underline{\sigma}, \bar{\sigma})>0 $ such that  
		\begin{equation*}
			\left|Y_{t}^{l}\right| \leq \exp \left(C(\bar{\sigma}, \underline{\sigma}%
			)\left(1+\Lambda^{T}(u, v)\right)\right)\left(\hat{\mathbb{E}}_{t}\left[%
			|\xi|^{2}\right]^{\frac{1}{2}}+\hat{\mathbb{E}}_{t}\left[\left(\int_{t}^{T}
			h_{s}^{l, 0} d s\right)^{2}\right]^{\frac{1}{2}}\right),  
		\end{equation*}
		where $h_s^{l,0} =\left|f^l(s,0,0)\right|+$ $\left|g^l(s,0,0)\right|$. 
	\end{proposition}
	
	\begin{theorem}
		Given that $f^l $, $g^l $
		satisfy the assumptions $ \textnormal{(A2)-(A3)} $ and $\xi^l \in L_{G}^{2+\beta}\left(\Omega_{T}\right) $ for some $\beta>0 $. And suppose that $\left(Y^{l}, Z^{l}, K^{l}\right) \in \mathfrak{S}_{G}^{2}(0, T), l=1,2 $
		is the solution for the following $G $-BSDEs on time interval $[0, T] $
		:  
		\begin{equation*}
			Y_{t}^{l}=\xi^{l}+\int_{t}^{T} f^{l}\left(s, Y_{s}^{l}, Z_{s}^{l}\right) d
			s+\sum_{i, j=1}^{d} \int_{t}^{T} g_{i j}^{l}\left(s, Y_{s}^{l},
			Z_{s}^{l}\right) d\left\langle B^{i}, B^{j}\right\rangle_{s}-\int_{t}^{T}
			Z_{s}^{l} d B_{s}-\left(K_{T}^{l}-K_{t}^{l}\right).
		\end{equation*}
		If $\xi^{1} \geq \xi^{2} $ and  
		\begin{equation*}
			f^{1}\left(s, Y_{s}^{2}, Z_{s}^{2}\right) \geq f^{2}\left(s, Y_{s}^{2},
			Z_{s}^{2}\right),\left(\tilde{g}_{i j}^{1}\left(s, Y_{s}^{2},
			Z_{s}^{2}\right)\right)_{i, j=1}^{d} \geq\left(\tilde{g}_{i j}^{2}\left(s,
			Y_{s}^{2}, Z_{s}^{2}\right)\right)_{i, j=1}^{d} \text { with } \tilde{g}_{i
				j}^{l}=\frac{g_{i j}^{l}+g_{j i}^{l}}{2} \text {, } 
		\end{equation*}
		then $Y_{t}^{1} \geq Y_{t}^{2}$. 
	\end{theorem}
	
	The explicit solution of linear $G $-BSDE has always been a crucial
	foundational issue, playing an important role in the study of the solution of  $G $-BSDE. Therefore, we will focus on the linear $G $-BSDE with time-varying Lipschitz generators. Consider the following linear $G $-BSDE  
	\begin{equation}  \label{eq9}
		Y_{t}=\xi +\int_{t}^{T} f(s)ds+\int_{t}^{T} g(s)d\langle B \rangle_{s}
		-\int_{t}^{T}Z_{s}dB_{s}-(K_{T}-K_{t}),
	\end{equation}
	where the generators satisfy the conditions
	\begin{equation*}
		f(s)=a_sY_s+b_sZ_s+m_s ,\text{ } g(s)=a_sY_s+b_sZ_s+m_s, 
	\end{equation*}
	with $a_s$, $c_s$, $m_s$, $n_s\in M_G^1(0,T)$, $b_s$, $d_s\in M_G^2(0,T)$, $\xi
	$, $\int_{0}^{T}|m_s|\mathrm{d}s$, $\int_{0}^{T}|n_s|\mathrm{d}s\in
	L_G^2(\Omega_{T})$. Additionally, assume that $a_s$, $b_s$, $c_s$, $d_s$ satisfy the
	assumption (A2) and (A3).
	
	The linear $ G $-BSDE is different from the classical case, the explicit solution of it under $G $-expection can only be provided through an auxiliary $\widetilde{G} $-expection space $(\widetilde{\Omega}, L_{\widetilde{G%
	}}^{1}(\widetilde{\Omega}), \hat{\mathbb{E}}^{\widetilde{G}})$. More
	detailed information can be found in \cite{hu2020bsdes}. Set $\widetilde{\Omega}=C_{0}([0, \infty), \mathbb{R}^{2 d})$ and  
	\[
	\widetilde{G}(A)=\frac{1}{2} \sup _{Q \in \Gamma} \operatorname{tr}\left[A\left[\renewcommand{\arraystretch}{0.8} \begin{array}{cc} 
		Q & I_{d} \\ 
		I_{d} & Q^{-1} 
	\end{array}\right]\right], \quad \forall A \in \mathbb{S}_{2 d}, \text{ } \Gamma \subset \mathbb{S}_{d}^+,
	\]
	where $\Gamma$ is defined in \eqref{eq26} and the process $	(B_t,\widetilde{B}_{t})_{t \geq 0}$ is a canonical process on the extended $\widetilde{G} $-expectation space $(\widetilde{\Omega}, L_{\widetilde{G}}^{1}(\widetilde{\Omega%
	}), \hat{\mathbb{E}}^{\widetilde{G}})$. And let $\Gamma_{t} $ be the solution of the following linear $\widetilde{G} $-SDE:
	\begin{equation}\label{eq27}
	\Gamma_{t}=1+\int_{0}^{t} a_{s} \Gamma_{s} d s+ \int_{0}^{t} c_{s}
	\Gamma_{s} d\left\langle B, B\right\rangle_{s}+\int_{0}^{t} d_{s} \Gamma_{s}
	d B_{s}+\int_{0}^{t} b_{s} \Gamma_{s} d \widetilde{B}_{s},
	\end{equation}
    which asists in obtaining the solution of
    the linear $G $-BSDE $\eqref{eq9}$. The relevant results are as follows.
	\begin{proposition}
		Suppose that the assumptions $ \textnormal{(A2)-(A3)} $ are satisfied. And let $(Y, Z, K) \in \mathfrak{S}%
		_{G}^{2}(0, T) $ be the solution of	linear $G$-BSDE \eqref{eq9} on [0,T]. Then, in the extended $\widetilde{G} $-expectation space,
		we can get  
		\begin{equation*}
			Y_{t}=\Gamma_{t}^{-1} \hat{\mathbb{E}}_{t}^{\tilde{G}}\left[\Gamma_{T}
			\xi+\int_{t}^{T} m_{s} \Gamma_{s} d s+\sum_{i, j=1}^{d} \int_{t}^{T}
			n_{s}^{i j} \Gamma_{s} d\left\langle B^{i}, B^{j}\right\rangle_{s}\right],  
		\end{equation*}
		where $\left\{\Gamma_{t}\right\}_{t \in[0, T]} $ is the solution of the
		linear $\widetilde{G} $-SDE \eqref{eq27}. Moreover,  
		\begin{equation*}
			\Gamma_{t}^{-1} \hat{\mathbb{E}}_{t}^{\widetilde{G}}\left[\Gamma_{T}
			K_{T}-\int_{t}^{T} a_{s} K_{s} \Gamma_{s} d s-\sum_{i, j=1}^{d} \int_{t}^{T}
			c_{s}^{i j} K_{s} \Gamma_{s} d\left\langle B^{i}, B^{j}\right\rangle_{s}%
			\right]=K_{t} .  
		\end{equation*}
	\end{proposition}
	\begin{proposition}
		\textnormal{\cite{wang2021backward}} Assume that $ X_n $, $ n\geq 1 $ and $ X $ are $ \mathcal{B}(\Omega) $-mesurable. Let $ \left\{X_n\right\}_{n=1}^\infty$ be such that $ X_n\downarrow X$, q.s. Then $ \hat{\mathbb{E}}[X_n ] \downarrow \hat{\mathbb{E}}[X]. $
	\end{proposition}
	\section{The Well-posedness of \texorpdfstring{$ G $}{}-BSDEs with Time-varying Uniformly continuous generators}
	
	In this section, we consider the well-posedness of following one-dimensional $G $-BSDE \eqref{eq1} and the corresponding comparison theorem on the finite interval [0, T]: 
	\begin{equation}\label{eq1} 
		Y_{t}=\xi +\int_{t}^{T} f(s,Y_{s},Z_{s})ds+\int_{t}^{T}
		g(s,Y_{s},Z_{s})d\langle B \rangle_{s}
		-\int_{t}^{T}Z_{s}dB_{s}-(K_{T}-K_{t}),
	\end{equation}
	where the generators  
	\begin{equation*}
		f(t, \omega, y, z),\text{ }g(t, \omega, y, z):[0, T] \times \Omega \times \mathbb{R}
		\times \mathbb{R}^{d} \rightarrow \mathbb{R}  
	\end{equation*}
	and the terminal $\xi $ satisfy the following assumptions.
	\vspace{1em}  \newline
	(H1)  There exists a constant $\beta > 2$ such that for each $\left(y, z\right)\in \mathbb{R} \times \mathbb{R}^{d}$, $f(\cdot ,\cdot,y,z)$, $g(\cdot ,\cdot ,y,z)\in M_{G}^{1}(0,T) $ and  
	\begin{equation*}
		\int_{0}^{T}(|f(s,0,0)|+|g(s,0,0)|)ds\in L_{G}^{\beta } (\Omega _{T}).
	\end{equation*}
	(H2) There are two positive non-random processes $u(\cdot)$ and $v(\cdot)$: $[0,T] \rightarrow \mathbb{R}^{+}$, and a continuous function $\phi(\cdot)$ with a constant $L$, independent of $(t, \omega)$, such that 
	\begin{equation*}
		|f(t,\omega,y_{1},z_{1})-f(t,\omega,y_{2},z_{2})|+|g(t,\omega,y_{1},z_{1})-g(t,\omega,y_{2},z_{2})|\le u(t)|y_{1}-y_{2}|+v(t)\phi (|z_{1}-z_{2}|),  
	\end{equation*}
	where $\phi (\cdot ) :\mathbb{R}^{+} \rightarrow \mathbb{R}^{+} $ is a
	non-decreasing and sub-additive function satisfying that $\phi(z) \leq
	L(1+z) $ as well as $\phi(0) = 0 $.\newline
	(H3) The non-random process $u(t)$, $v(t)$ satisfy the condition of integrability :  
	\begin{equation*}
		\wedge (u,v)=\int_{0}^{T} (u(t)+v^{2}(t))dt< \infty.  
	\end{equation*}
	\begin{remark}
		\textnormal{(1)} The condition \textnormal{(H2)} can be converted to  
		\begin{equation*}
			|f(t,\omega,y_{1},z_{1})-f(t,\omega,y_{2},z_{2})|+|g(t,\omega,y_{1},z_{1})-g(t,\omega,y_{2},z_{2})|  \leq u(t)|y_{1}-y_{2}|+Lv(t)(1+|z_{1}-z_{2}|)  
		\end{equation*}
		For the BSDE case with uniformly continuous generators, more details can be seen in \textnormal{\cite{mao1995adapted}}.
		
		\textnormal{(2)} To ensure the convergence, for each $t \in [0,T] $, $u(t)$ and $v(t) $ need to be finite function. 
	\end{remark}
	
	In order to prove the existence and uniqueness of the solution to $G $-BSDE \eqref{eq1}, we need to prove some useful Lemmas as follows. Next, we first construct sequences of Lipschitz function $\underline{\varphi}_{n}(t,\omega ,y,z)$ and $\bar{\varphi}_{n}(t,\omega ,y,z)$ which can
	approximates the generators $f$ and $g$ respectively. And based on them, the supremum approximation equations and the infimum approximation equations of $G$-BSDE \eqref{eq1} are derived (see \eqref{eq5} \eqref{eq6}). Therefore, for any $(t,\omega,y,z) \in [0,T]\times \Omega \times \mathbb{R} \times \mathbb{R}^{d} $, any $n\in N^{+}$, denote  
	\begin{equation}  \label{eq_xia}
		\underline{\varphi}_{n}(t, \omega, y, z):=\inf _{q\in Q}\left.\{\varphi(t,
		\omega, y,q)+nv(t)|z-q|\right\}-\varphi_{0}(t),
	\end{equation}
	\begin{equation}  \label{eq_shang}
		\overline{\varphi}_{n}(t, \omega, y, z):=\sup _{q \in Q}\left.\{\varphi(t,
		\omega, y, q)-nv(t)|z-q|\right\}-\varphi_{0}(t),
	\end{equation}
	where $\varphi_{0}(t)=\varphi(t,0,0)$ and  $\varphi =f$, $g. $ The main properties
	of $\underline{\varphi}_{n} $ and $\overline{\varphi}_{n}$ are outlined as
	follows.
	\begin{lemma}
		Set $L=\max\left\{L,\text{ }1 \right\}, $ where $L $ is the Lipschitz constant associated with the function $\phi $. Under the assumptions \textnormal{(H1)-(H2)}, it can be inferred that for
		any $n>L$, the properties \textnormal{(\romannumeral 1)}-\textnormal{(\romannumeral 6)} hold:\\
		\textnormal{(\romannumeral 1)} $\underline{\varphi}_{n} $ and $\overline{\varphi}_{n} $ satisfy the linear growth
		condition, that is, for any $(t, \omega, y, z)$, 
		\begin{align*}
			-L&(u(t)|y|+v(t)|z|+v(t)) \leq \underline{\varphi}_{n}(t, \omega, y, z) \leq
			\varphi(t, \omega, y, z) \\
			&-\varphi_{0}(t)\leq \overline{\varphi}_{n}(t, \omega, y, z) \leq
			L(u(t)|y|+v(t)|z|+v(t)),
		\end{align*}
		where $\varphi =f$, $ g$. \newline
		\textnormal{(\romannumeral 2)} For all $(t, \omega, y, z)$, the function $\left\{\underline{\varphi}%
		_{n}(t, \omega, y, z)\right\}_{n \geq 1}$ is non-decreasing with respect to $%
		n$ and $\left\{\bar{\varphi}_{n}(t, \omega, y, z)\vphantom{\underline{%
				\varphi}_{n}(t, \omega, y, z)}\right\}_{n \geq 1}$ is non-increasing with
		respect to $n$.\newline
		\textnormal{(\romannumeral 3)}The sequence of functions $\left\{\underline{\varphi}_{n}(t, \omega, y,
		z)\right\}_{n \geq 1}$ and $\left\{\overline{\varphi}_{n}(t, \omega, y, z)%
		\vphantom{\underline{\varphi}_{n}(t, \omega, y, z)}\right\}_{n \geq 1}$ with
		respect to the variables $(y, z) $ are Lipschitz functions. Then, for any $t \in [0, T]$, $\omega \in \Omega_{T}$, $y_{1}$, $
		y_{2} \in \mathbb{R}$, $z_{1}$, $ z_{2} \in \mathbb{R}^{d},$ 
		\begin{equation*}
			\left|\underline{\varphi}_{n}\left(t,\omega,y_{1},z_{1}\right)-\underline{%
				\varphi}_{n}\left(t,\omega,y_{2}, z_{2}\right)\right| \leq
			u(t)\left|y_{1}-y_{2}\right|+n v(t)\left|z_{1}-z_{2}\right|,
		\end{equation*}
		\begin{equation*}
			\left|\bar{\varphi}_{n}\left(\omega, t, y_{1}, z_{1}\right)-\bar{\varphi}%
			_{n}\left(\omega, t, y_{2}, z_{2}\right)\right| \leq
			u(t)\left|y_{1}-y_{2}\right|+n v(t)\left|z_{1}-z_{2}\right|. 
		\end{equation*}
		\textnormal{(\romannumeral 4)} $\underline{\varphi}_{n} $ and $\overline{\varphi}_{n} $ satisfy the
		time-varying uniformly continuous condition about $z$, i.e. for any $t \in
		[0,T]$, $\omega \in \Omega_{T}$, $y\in \mathbb{R}$, $z_{1}$, $z_{2}\in \mathbb{R}^{d},$ 
		\begin{equation*}
			\left|\underline{\varphi}_{n}\left(t,\omega, y, z_{1}\right)-\underline{%
				\varphi}_{n}\left(t,\omega, y, z_{2}\right)\right| \leq v(t)
			\phi\left(\left|z_{1}-z_{2}\right|\right), 
		\end{equation*}
		\begin{equation*}
			\left|\bar{\varphi}_{n}\left(\omega, t, y, z_{1}\right)-\bar{\varphi}%
			_{n}\left(\omega, t, y, z_{2}\right)\right| \leq v(t)
			\phi\left(\left|z_{1}-z_{2}\right|\right). 
		\end{equation*}
		\textnormal{(\romannumeral 5)} For all $(t, \omega, y, z)\in [0, T] \times \Omega \times \mathbb{R}
		\times \mathbb{R}^{d} $, 
		\begin{equation*}
			0 \leq \varphi(t,\omega, y, z)-\varphi_{0}(t)-\underline{\varphi}%
			_{n}(t,\omega, y, z) \leq v(t) \phi \left(\frac{2 L}{n-L}\right), 
		\end{equation*}
		\begin{equation*}
			0 \leq \bar{\varphi}_{n}(t,\omega, y, z)-\varphi(t,\omega, y,
			z)+\varphi_{0}(t) \leq v(t) \phi \left(\frac{2 L}{n-L}\right). 
		\end{equation*}
		\textnormal{(\romannumeral 6)} If $\left(y_ {n}, z_ {n}\right) \rightarrow(y, z)$ as $n\rightarrow\infty $, then 
		\begin{align*}
			&\underline{\varphi}_{n}\left(t, \omega, y_ {n}, z_{n}\right) \rightarrow
			\varphi(t, \omega, y, z)-\varphi_{0}(t), \\
			&\bar{\varphi}_{n}\left(t, \omega, y_{n}, z_{n}\right) \rightarrow
			\varphi(t, \omega, y, z)-\varphi_{0}(t).
		\end{align*}
		\begin{proof}
			\textnormal{(\romannumeral 1)} Based on \textnormal{(H2)}, it is obvious that for any $n>L\geq1$,
			\begin{align*}
				\bar{\varphi}_{n}(t,\omega ,y,z)
				&\leq\sup _{q \in Q}\left\{\varphi(t,\omega,y,q)-\varphi(t,0,0)-L v(t)|z-q|\right\}\\
				&\leq \sup _{q \in Q}\left\{Lu(t)|y|+Lv(t)(1+|q|)-Lv(t)|z-q|\right\}\\
				&\leq L(u(t)|y|+v(t)|z|+v(t))
			\end{align*}
			and 
			\begin{align*}
				\underline{\varphi}_{n}(t, \omega, y, z)&\geq\inf _{q \in Q}\left\{\varphi(t, \omega, y, q)-\varphi(t,0,0)+L v(t)|z-q|\right\}\\
				&\ge \inf _{q\in Q}\left\{-Lu(t)|y|-Lv(t)(1+|q|)+Lv(t)|z-q|\right\}\\
				&\ge -L\left(u(t)|y|+v(t)|z|+v(t)\right).
			\end{align*}
			Note that for each fixed $ (n,y,z) $, each $ q\in Q $, we have
			$
			\varphi(t,\omega,y,q)\pm nv(t)|z-q|\rightarrow \varphi(t,\omega,y,z)
			$ as $ q\rightarrow z $. Then for any $ \varepsilon >0$, there exist $ q_1 $, $ q_2\in Q $ such that
			\[
			\inf _{q \in Q}\left\{\varphi(t, \omega, y, q)+nv(t)|z-q|\right\}-\varepsilon
			\leq \varphi(t, \omega, y, q_1)+nv(t)|z-q_1|-\varepsilon
			\leq \varphi(t, \omega, y, z),
			\]
			\[
			\sup _{q \in Q}\left\{\varphi(t, \omega, y, q)-nv(t)|z-q|\right\}+\varepsilon
			\geq \varphi(t, \omega, y, q_2)-nv(t)|z-q_2|+\varepsilon
			\geq \varphi(t, \omega, y, z).
			\]
			Therefore, sending $ \varepsilon \rightarrow 0 $, we can obtain \[ \underline{\varphi}_{n}(t, \omega, y, z) \leq
			\varphi(t, \omega, y, z)-\varphi_{0}(t)\leq \overline{\varphi}_{n}(t, \omega, y, z), \] which  ends the proof.
			\\
			\textnormal{(\romannumeral 4)} According to the definitions of $\underline{\varphi}_{n}(t,\omega ,y,z)$, it can be inferred that $ \underline{\varphi}_{n} $ can be converted to
			\begin{align*}
				\underline{\varphi}_{n}(t, \omega, y, z) =\inf _{q\in Q}\{\varphi(t, \omega, y,z-q)+nv(t)|q|\}-\varphi_{0}(t). 
			\end{align*}
			Recalling inequality \[|\inf_{x \in D}\left\{f_{1}(x)\right\}-\inf_{x \in D}\left\{f_{2}(x)\right\}|\le \sup _{x \in D}\left\{\left|f_{1}(x)-f_{2}(x)\right|\right\},\] we can get that
			\begin{align*}
				\left|\underline{\varphi}_{n}\left(\omega, t, y, z_{1}\right)-\underline{\varphi}_{n}\left(\omega, t, y, z_{2}\right)\right|
				&\le 
				\sup _{ q \in \mathbf{R} }\{|\varphi(t, \omega, y, z_{1}-q)-\varphi(t, \omega, y, z_{2}-q)|\}\\
				&\leq 
				v(t) \phi\left(\left|z_{1}-z_{2}\right|\right).
			\end{align*}
			Therefore, the first formula has been proven. And the second one can be proved similarly.\\  
			\textnormal{(\romannumeral 5)} To simplify, we will only prove the first formula. For all $ (t, \omega, y, z) $,
			\begin{equation}\label{eq4}
				\varphi(t,\omega, y, q) \\
				\geq \varphi(t,\omega, y, z)-v(t)\phi(|z-q|) \\
				\geq \varphi(t,\omega, y, z)-Lv(t)(1+|z-q|), \quad \forall q\in Q .
			\end{equation}
			For any $n>L$, define two sets:
			\[
			\Lambda_{n}:= \left\{ q \in Q : n \lvert z-q \rvert \geq L(\lvert z-q \rvert + 2) \right\},
			\]
			\[
			\Lambda_{n}^{c}:=\left\{q\in Q: n\lvert z-q \rvert < L(\lvert z-q \rvert +2) \right\}.
			\]
			It is easy to note that the sets $\Lambda_{n}$ and $\Lambda_{n}^{c}$ are both non-empty and $Q = \Lambda_{n} \cup \Lambda_{n}^{c}$. Suppose $q\in \Lambda_{n}$, in accordance with the definition of $\Lambda_{n}$ and inequality \eqref{eq4}, we conclude that
			\begin{align*}
				\varphi(t,\omega,y,q)-\varphi_{0}(t)+nv(t)|z-q|
				& \geq \varphi(t,\omega, y, q)-\varphi_{0}(t)+Lv(t)(|z-q|+2) \\
				& \geq \varphi(t,\omega, y, z)-L v(t)(|z-q|+1)-\varphi_{0}(t)+Lv(t)(|z-q|+2)\\
				& = \varphi(t,\omega, y, z)-\varphi_{0}(t)+L v(t)\\
				& >\underline{\varphi}_{n}(t,\omega, y, z) .
			\end{align*}
			Using the non-decreasing property of function $ \phi $ and the above inequality, it can be obtained that
			\begin{align*}
				\underline{\varphi}_{n}(t,\omega, y, z)
				& =\inf _{q \in \Lambda_{n} \cup \Lambda_{n}^{c}}\{\varphi(t,\omega, y, q)+n v(t)|z-q|\}-\varphi_{0}(t) \\
				& =\inf _{q \in \Lambda_{n}^{c}}\{\varphi(t,\omega, y, q)+n v(t)|z-q|\} -\varphi_{0}(t)\\
				& \geq \inf _{q \in \Lambda_{n}^{c}}\left\{\varphi(t,\omega, y, q): q \in Q ,|z-q|<\frac{2L}{n-L} \right\}-\varphi_{0}(t) \\
				& \geq \inf _{q\in \Lambda_{n}^{c} }\left\{\varphi(t,\omega, y, z)-v(t)\phi(|z-q|): q\in Q,|z-q|<\frac{2L}{n-L} \right\} -\varphi_{0}(t)\\
				& =\varphi(t,\omega, y, z)-v(t) \phi(\frac{2 L}{n-L})-\varphi_{0}(t). 
			\end{align*}
			Thus the proof of \textnormal{(\romannumeral 5)} is complete. And \textnormal{(\romannumeral 6)} can be derived from \textnormal{(\romannumeral 5)}.
		\end{proof}
	\end{lemma}
	
	Based on the generator $\underline{\varphi}_{n}$ and $\bar{\varphi}_{n} $, we construct the supremum approximation equation and the infimum
	approximation equation of $G $-BSDE \eqref{eq1}, which is defined as follows: 
	\begin{equation}  \label{eq5}
		\begin{aligned} \underline{Y}_{t}^{n}= &
			\xi+\int_{t}^{T}\left[\underline{f}_{n}\left(s, \underline{Y}_{s}^{n},
			\underline{Z}_{s}^{n}\right)+f_{0}(s)\right] d
			s+\int_{t}^{T}\left[\underline{g}_{n}\left(s, \underline{Y}_{s}^{n},
			\underline{Z}_{s}^{n}\right)+g(s)\right] d\left\langle
			B\right\rangle_{s} \\& -\int_{t}^{T} \underline{Z}_{s}^{n} d
			B_{s}-\left(\underline{K}_{T}^{n}-\underline{K}_{t}^{n}\right), \end{aligned}
	\end{equation}
	\begin{equation}\label{eq6}
		\begin{aligned} \bar{Y}_{t}^{n}= & \xi+\int_{t}^{T}\left[\bar{f}_{n}\left(s,
			\bar{Y}_{s}^{n}, \bar{Z}_{s}^{n}\right)+f_{0}(s)\right] d
			s+\int_{t}^{T}\left[\bar{g}_{n}\left(s, \bar{Y}_{s}^{n},
			\bar{Z}_{s}^{n}\right)+g_{0}(s)\right] d\left\langle
			B\right\rangle_{s} \\& -\int_{t}^{T} \bar{Z}_{s}^{n} d
			B_{s}-\left(\bar{K}_{T}^{n}-\bar{K}_{t}^{n}\right). \end{aligned}
	\end{equation}
	
	To ensure the existence of $\underline{Y}^{n}$ and $\bar{Y}^{n}$, assumption (H4) is proposed. Following this, some
	properties and a priori estimates of the $G $-BSDEs (                                    $\underline{Y}_{t}^{n}$
	and $\bar{Y}_{t}^{n}$) are given.  \vspace{1em}  \newline
	(H4) For each $n
	\in \mathbb{N} $ and for each $(y,z) \in \mathbb{R} \times \mathbb{R}^{d} $, $\underline{\varphi}_{n}(\cdot, \cdot, y, z)$, $\bar{\varphi}_{n}(\cdot, \cdot, y, z)\in M_{G}^{1 }(0,T)$ where $\varphi=f$, $g$.
	\begin{remark}
		\textnormal{(1)} The assumption \textnormal{(H4)} can be proven to hold in many situations. For instance, in the case where the functions $ f $ and $ g $ are uniformly continuous with respect to $ (t, w) $, since
		\begin{align*}
			\underline{\varphi}_{n}(t_1, \omega_1, y, z)-\underline{\varphi}_{n}(t_2, \omega_2, y, z)
			&\leq|\inf _{q \in Q}\left\{\varphi(t_1, \omega_1, y, q)+nv(t)|z-q|\right\}-\inf _{q \in Q}\left\{\varphi(t_2, \omega_2, y, q)+nv(t)|z-q|\right\}|\\
			&\le \sup_{q \in Q}\left\{\varphi(t_1, \omega_1, y, q)-\varphi(t_2, \omega_2, y, q)\right\}, \quad \varphi=f,g,
		\end{align*}
		we can directly get $ \underline{\varphi}_{n} $ and $ \bar{\varphi}_{n} $ are also uniformly continuous with respect to $ (t, w) $. Following this, by \textnormal{(\romannumeral 1)} of Lemma \textnormal{3.2} and Theorem \textnormal{4.7} in \textnormal{\cite{hu2016quasi}}, we have $\underline{\varphi}_{n}(\cdot, \cdot, y, z)$, $\bar{\varphi}_{n}(\cdot, \cdot, y, z)\in M_{G}^{1 }(0,T)$ for each $(y,z)$. Then the assumtion \textnormal{(H4)} can be easily verified.\\
		\textnormal{(2)} The $C(\underline{\sigma},\bar{\sigma},\Lambda (u,v))$ in following Lemmas is a constant depending on $\underline{\sigma},\bar{\sigma}$ and $\Lambda (u,v)$, which may vary from line to line.
	\end{remark}  
	\begin{lemma}
		Suppose that $\xi \in L_{G}^{\beta}(\Omega _{T})$ for some $\beta > 2 $ and
		the assumptions \textnormal{(H1)}-\textnormal{(H4)} are satisfied. Then G-BSDE \eqref{eq5}$ (resp.\eqref{eq6}) $ admits a unique solution $(\bar{Y}, \bar{Z}, \bar{K})$$ (resp. (\underline{Y}, \underline{Z}, \underline{K}) ) $  in $\mathfrak{S}_{G}^{2}(0, T)$. Moreover, the inequalities $\underline{Y}^{n} \leq \underline{Y}^{n+1} \leq \bar{Y}^{m+1} \leq \bar{Y}^{m}$ hold for each $n$, $ m \in \mathbb{N}$ and $\underline{Y}^{n}$, $\bar{Y}^{n}$ are uniformly bounded with respect to n in $S_{G}^{2}(0, T)$.
		
		\begin{proof}
			(1)Recalling Theorem \textnormal{2.3} and Lemma \textnormal{3.2}, it can be deduced that there exists a unique solution for the G-BSDE \eqref{eq5} and G-BSDE \eqref{eq6} respectively. Meanwhile, we conclude  that 
			\[
			\underline{f}^n \leq \underline{f}^{n+1} \leq \bar{f}^m \leq \bar{f}^{m+1}\hspace{0.3cm}\text{and}\hspace{0.3cm}\underline{g}^n \leq \underline{g}^{n+1} \leq \bar{g}^m \leq \bar{g}^{m+1}.
			\]Furthermore, according to the comparison theorem( see Theorem \textnormal{2.5}), we have \[
			\underline{Y}^{n} \leq \underline{Y}^{n+1} \leq \bar{Y}^{m+1} \leq \bar{Y}^{m},\quad \forall n, m \in \mathbb{N}.
			\]
			(2)To obtain a uniform boundedness estimate for $\underline{Y}^{n}$ and $\bar{Y}^{n}$, we construct the following $ G $-BSDE. Let $R(t, y, z)=L\left(u(t)|y|+v(t)|z|+v(t)\right)$,   
			\begin{equation}\label{eq7}
				\begin{aligned}
					U_{t}= & \xi+\int_{t}^{T}\left[R\left(U_{s}, V_{s}\right)+f_{0}(s)\right] d s+\int_{t}^{T}\left[R\left(U_{s}, V_{s}\right)+g_{0}(s)\right] d\left\langle B\right\rangle_{s} \\& -\int_{t}^{T} V_{s} d B_{s}-\left(R_{T}-R_{t}\right),
					\\U_{t}^{\prime}= & \xi+\int_{t}^{T}\left[-R\left(U_{s}^{\prime}, V_{s}^{\prime}\right)+f_{0}(s)\right] d s+\int_{t}^{T}\left[-R\left(U_{s}^{\prime}, V_{s}^{\prime}\right)+g_{0}(s)\right] d\left\langle B\right\rangle_{s} \\& -\int_{t}^{T} V_{s}^{\prime} d B_{s}-\left(R_{T}^{\prime}-R_{t}^{\prime}\right) .
				\end{aligned}
			\end{equation}
			Similar to the proof of (1), we can obtain $ R(t, y, z) \in M_G^1(0,T) $ and G-BSDE \eqref{eq7} has a unique solution $(U,V,R)$, $(U^{\prime},V^{\prime},R^{\prime})\in \mathfrak{S}_{G}^{2}(0, T)$.  Moreover, it holds that for any $ n\in \mathbb{N} $,
			\[
			U_{t}^{\prime}\leq \underline{Y}^{n}_{t} \leq \bar{Y}^{n}_{t} \leq U_{t}, \quad \forall t \in [0,T].
			\]
			According to Proposition \textnormal{2.4}, we can get that the uniform boundedness of $\underline{Y}^{n}$ and $\bar{Y}^{n}$ in $S_{G}^{2}(0, T)$. The proof is complete.
		\end{proof}
	\end{lemma}
	
	The following Lemma is a fundamental result of this paper. The main idea is to obtain an estimate of $\bar{Y}^{n} - \underline{Y}^{n}$ through the linearization of $G$-BSDE, which can assist in proving the subsequent convergence. In this proof, we construct a Lipschitz continuous linear decomposition function $ l(y) $, which is suitable for the time-varying condition.  The function $ l(y) $ plays a crucial role in the subsequent proof.	
	\begin{lemma}
		For each $n>L$ and each $t \in [0, T]$, the difference between $\underline{Y}^{n}$
		and $\bar{Y}^{n}$ can be uniformly controlled, as the following
		inequality holds:  
		\begin{equation*}
			\left|\bar{Y}_{t}^{n} - \underline{Y}_{t}^{n}\right| \leq C\left(\underline{\sigma},\bar{\sigma},\Lambda (u,v)\right)\phi\left(\frac{2 L}{n-L}\right),\quad
			\forall t \in[0, T],  
		\end{equation*}
		where  $C(\underline{\sigma},\bar{\sigma},\Lambda (u,v))$ is a constant depending on $\underline{\sigma},\bar{\sigma}$ and $\Lambda (u,v)$.  
		
		\begin{proof}
			Set $(\hat{Y}, \hat{Z})=\left(\bar{Y}^{n}-\underline{Y}^{n}, \bar{Z}^{n}-\underline{Z}^{n}\right) $. For each $t \in [0, T]$, we have
			\\
			\begin{equation}\label{eq8}   \hat{Y}_t+\underline{K}_t^n=\underline{K}_T^n+\int_t^{T} \hat{f}_s ds+\int_t^{T} \hat{\mathrm{g}}_s d\langle B\rangle_s-\int_t^{T} \hat{Z}_s d B_s-\left(\bar{K}_T^n-\bar{K}_t^n\right),
			\end{equation}
			\leavevmode \\
			where $\hat{\varphi}_s=\bar{\varphi}_n\left(s, \bar{Y}_s^n, \bar{Z}_s^n\right)-\underline{\varphi}_n\left(s, \underline{Y}_s^n, \underline{Z}_s^n\right)$ and $\varphi=f$, $ g$. Next we linearize the above $G$-BSDE \eqref{eq8} in the following way. Set $ m_{s}=\bar{f}_n\left(s, \underline{Y}_s^n, \underline{Z}_s^n\right)-\underline{f}_n\left(s, \underline{Y}_s^n, \underline{Z}_s^n\right)$, $n_{s}=
			\bar{g}_n\left(s, \underline{Y}_s^n, \underline{Z}_s^n\right)-\underline{g}_n\left(s, \underline{Y}_s^n, \underline{Z}_s^n\right) $ and define function
			\[
			l(y) = \mathbf{I}_{\left\{|y| \leq \varepsilon\right\}} + \left(2 - \frac{|y|}{\varepsilon}\right) \mathbf{I}_{\left\{\varepsilon < |y| < 2\varepsilon \right\}}, \quad \forall y \in \mathbb{R}.
			\]
			Note that $ (1-l(y))y^{-1} $ is a bounded and Lipschitz continuous function with respect to $ y $. Then for any $ \varepsilon > 0$, there exist four processes $a^{\varepsilon }_s$, $ b^{\varepsilon }_s$, $c^{\varepsilon}_s$, $ d^{\varepsilon}_s $ such that for any $s \in[0, T]$, we can get that
			\[
			\hat{f}_s=a_s^{\varepsilon } \hat{Y}_s+b_s^{\varepsilon } \hat{Z}_s+m_s+m_s^{\varepsilon } ,\text{ }\hat{g}_s=c_s^{\varepsilon }  \hat{Y}_s+d_s^{\varepsilon } \hat{Z}_s+n_s+n_s^{\varepsilon },
			\]
			where 
			\[
			a_{s}^{\varepsilon}(\hat{Y}_s, \bar{Z}_s^n) = (1-l(\hat{Y}_s))\frac{\bar{f}_n(s, \bar{Y}_s^n, \bar{Z}_s^n) - \bar{f}_n(s, \underline{Y}_s^n, \bar{Z}_s^n)}{\hat{Y}_s}, \quad \text{if} \ \hat{Y}_s \neq 0,
			\]
			\[
			b_{s}^{\varepsilon}(\hat{Z}_s, \underline{Y}_s^n) = (1-l(\hat{Z}_s))\frac{\bar{f}_n(s, \underline{Y}_s^n, \bar{Z}_s^n) - \bar{f}_n(s, \underline{Y}_s^n, \underline{Z}_s^n)}{\hat{Z}_s}, \quad \text{if} \ \hat{Z}_s \neq 0.
			\]
			And when $ \hat{Y}_s $ and $ \hat{Z}_s $ are equal to 0, set $ a_{s}^{\varepsilon} = b_{s}^{\varepsilon}=0 $. It is evident that $a_{s}^{\varepsilon}(\hat{Y}_s, \bar{Z}_s^n)$ and $c_{s}^{\varepsilon}(\hat{Y}_s, \bar{Z}_s^n)$ are continuous functions with respect to $(\bar{Y}_s^n, \underline{Y}_s^n, \bar{Z}_s^n)$, and the same applies to $b^{\varepsilon }_s$, $d^{\varepsilon}_s$. Then we will show the bound of $a^{\varepsilon }_s$, $b^{\varepsilon }_s$, $m_s$,  $m^{\varepsilon}_s $. By \textnormal{(\romannumeral 3)} of Lemma \textnormal{3.2}, it is easy to verify that 
			$|a_{s}^{\varepsilon}| \leq u(s)$, $|b_{s}^{\varepsilon}| \leq nv(s)$. Note that $ m_s^{\varepsilon}= l(\hat{Y}_s)(\bar{f}_n(s, \bar{Y}_s^n, \bar{Z}_s^n) - \bar{f}_n(s, \underline{Y}_s^n, \bar{Z}_s^n))+l(\hat{Z}_s)(\bar{f}_n(s, \underline{Y}_s^n, \bar{Z}_s^n) - \bar{f}_n(s, \underline{Y}_s^n, \underline{Z}_s^n))$. It follows that
			\begin{equation*}
				\begin{aligned}
					|m_s^{\varepsilon }|&\leq
					u(s)|\hat{Y}_s|l(\hat{Y}_s)+nv(s)|\hat{Z}_s|l(\hat{Z}_s)\\
					&\leq 2\varepsilon(u(s)+nv(s)).
				\end{aligned}
			\end{equation*}
			Applying Lemma \textnormal{3.2 (\romannumeral 5)} to $m_s$, we get that
			$$
			\begin{aligned}
				|m_s|&=\left|\bar{f}_n\left(s, \underline{Y}_s^n, \underline{Z}_s^n\right)-f\left(s, \underline{Y}_s^n, \underline{Z}_s^n\right)+f_0(s)+f\left(s, \underline{Y}_s^n, \underline{Z}_s^n\right)-f_0(s)-\underline{f}_n\left(s, \underline{Y}_s^n, \underline{Z}_s^n\right)\right|\\
				& \leq 2v(s)\phi(\frac{2L}{n-L}).
			\end{aligned}
			$$
			And the same is true for the properties of $ c_s^{\varepsilon}$, $d_s^{\varepsilon }$, $n_s$ and $n_s^{\varepsilon} $. Set $ q(y)=(1-l(y))y^{-1}  $. Based on the fact that the function $ q(y) $ is bounded, we can conclude that  $ a_s^{\varepsilon}$, $b_s^{\varepsilon }$, $ m_s $, $ m_s^{\varepsilon}$, $ c_s^{\varepsilon}$, $d_s^{\varepsilon }$, $ n_s $ and  $ n_s^{\varepsilon}$ are Lipschitz functions with time-varying coefficients $ u(s) $ and $nv(s) $. Moreover, from Lemma \textnormal{3.2} in \textnormal{\cite{hu2020bsdes}} and the boundness of $ q(y) $, we have  $a_{s}^{\varepsilon}$, $c_{s}^{\varepsilon}$, $b_{s}^{\varepsilon}$, $d_{s}^{\varepsilon}\in M_{G}^1(0,T)$, which implies that they are quasi-continuous.	Meanwhile, note that \[
			\hat{\mathbb{E}}\left[\int_{0}^{T}|b_{s}^{\varepsilon}|^2\mathbf{1}_{\left\{|b_{s}^{\varepsilon}| \geq N\right\}}ds \right]\leq \int_{0}^{T}n^2v^2(s)\mathbf{I}_{\left\{nv(s) \geq N\right\}}ds,\quad \forall N\in\mathbb{N}.
			\]
			Then, by applying Theorem \textnormal{4.7} in \textnormal{\cite{hu2016quasi}} and the monotone convergence theorem, we can derived that $b_{s}^{\varepsilon}$, $d_{s}^{\varepsilon}\in M_{G}^2(0,T)$.
			Following  Lemma \textnormal{3.2} in \textnormal{\cite{hu2020bsdes}} and the assumption \textnormal{(H4)}, we can also  get that $m_{s}$, $m_{s}^{\varepsilon}$, $n_{s}$, $n_{s}^{\varepsilon}\in M_{G}^1(0,T)$ and $\int_{0}^{T}|m_{s}|d s$, $\int_{0}^{T}|m_{s}^{\varepsilon}|d s$, $\int_{0}^{T}|n_{s}|d s$, $\int_{0}^{T}|n_{s}^{\varepsilon}|ds\in L_{G}^2(\Omega_{T})$.
			For each $ t \in [0,T] $, $ s\in [t,T] $, according to Proposition \textnormal{2.6} and simple calculation, we  conclude that
			\begin{equation}\label{eq18}
				\begin{aligned}
					\hat{Y}_t+\underline{K}_t^n= & (\tilde{\Gamma}_t^{\varepsilon})^{-1}\hat{\mathbb{E}}_t^{\widetilde{G}}\left[\tilde{\Gamma}_T^{\varepsilon} \underline{K}_T^n+\int_t^T\left(m_s+m_s^{\varepsilon}-a_s^{\varepsilon} \underline{K}_s^n\right) \tilde{\Gamma}_s^{\varepsilon} d s\right. \left.+\int_t^T\left(n_s+n_s^{\varepsilon}-c_s^{\varepsilon} \underline{K}_s^n\right) \tilde{\Gamma}_s^{\varepsilon} d\langle B\rangle_s\right]\\
					\leq& (\tilde{\Gamma}_t^{\varepsilon})^{-1}\hat{\mathbb{E}}_t^{\widetilde{G}}\left[\int_t^T\left(m_s+m_s^{\varepsilon}\right) \tilde{\Gamma}_s^{ \varepsilon} ds\right. \left.+\int_t^T\left(n_s+n_s^{\varepsilon}\right) \tilde{\Gamma}_s^{\varepsilon} d\langle B\rangle_s\right]\\
					&+(\tilde{\Gamma}_t^{\varepsilon})^{-1}\hat{\mathbb{E}}_t^{\widetilde{G}}\left[\tilde{\Gamma}_T^{\varepsilon} \underline{K}_T^n-\int_t^T\left(a_s^{\varepsilon} \underline{K}_s^n\right) \tilde{\Gamma}_s^{\varepsilon} d s\right. \left.-\int_t^T\left(c_s^{\varepsilon} \underline{K}_s^n\right) \tilde{\Gamma}_s^{\varepsilon} d\langle B\rangle_s\right]\\
					=&(\tilde{\Gamma}_t^{\varepsilon})^{-1}\hat{\mathbb{E}}_t^{\widetilde{G}}\left[\int_t^T\left(m_s+m_s^{\varepsilon}\right) \tilde{\Gamma}_s^{\varepsilon} ds\right. \left.+\int_t^T\left(n_s+n_s^{\varepsilon}\right) \tilde{\Gamma}_s^{\varepsilon} d\langle B\rangle_s\right]+\underline{K}_t^n,
				\end{aligned}
			\end{equation}
			where $\left\{\tilde{\Gamma}_t^{\varepsilon}\right\}_{0\leq t \leq T}$ is the solution of $\widetilde{G}$-SDE
			\[
			\tilde{\Gamma}_t^{\varepsilon}=1+\int_{0}^{t} a_{r}^{\varepsilon}
			\tilde{\Gamma}_{r}^{\varepsilon} dr+\int_{0}^{t} c_{r}^{\varepsilon} \tilde{\Gamma}_{r}^{\varepsilon} d\left\langle B\right\rangle_{r}+ \int_{0}^{t} d_{r}^{\varepsilon}\tilde{\Gamma}_{r}^{\varepsilon}  d B_{r}+\int_{0}^{t} b_{r}^{\varepsilon} \tilde{\Gamma}_{r}^{\varepsilon}  d \widetilde{B}_{r}.
			\]
			Applying $  G \text {-Itô's formula } $ to $ \tilde{\Gamma}_t^{\varepsilon} $, we have
			\[
			\tilde{\Gamma}_t^{\varepsilon}=\exp\left(\int_{0}^{t}a_{r}^{\varepsilon}dr+\int_{0}^{t}c_{r}^{\varepsilon}d\left\langle B\right\rangle_{r}\right)\bar{\Gamma}_t^{\varepsilon},
			\]
			where \[ \bar{\Gamma}_t^{\varepsilon}=\exp\left(\int_{0}^{t}d_{r}^{\varepsilon}dB_r-\frac{1}{2}\int_{0}^{t}(d_{r}^{\varepsilon})^2d\langle B\rangle_{r}+\int_{0}^{t}b_{r}^{\varepsilon}d \widetilde{B}_{r}-\frac{1}{2}\int_{0}^{t}(b_{r}^{\varepsilon})^2d\langle \widetilde{B}\rangle_{r}-\int_{0}^{t}b_{r}^{\varepsilon}d_{r}^{\varepsilon}dr\right)\] is a symmetric $ G $-maringale. Combining with the bound of $a_{s}^{\varepsilon}  $ and $c_{s}^{\varepsilon}  $, it follows that for any $ s\in[t,T] $,
			\begin{equation}\label{eq19}
				\begin{aligned}
					(\tilde{\Gamma}_t^{\varepsilon})^{-1}\hat{\mathbb{E}}_t^{\widetilde{G}}[\tilde{\Gamma}_s^{\varepsilon}]
					&=\hat{\mathbb{E}}_t^{\widetilde{G}}\left[\exp\left(\int_{t}^{s}a_{r}^{\varepsilon}dr+\int_{t}^{s}c_{r}^{\varepsilon}d\left\langle B\right\rangle_{r}\right)(\bar{\Gamma}_t^{\varepsilon})^{-1}\bar{\Gamma}_s^{\varepsilon}\right]\\
					&\leq\hat{\mathbb{E}}_t^{\widetilde{G}}\left[\exp\left((1+\bar{\sigma}^2)\int_{0}^{T}u(s)ds\right)(\bar{\Gamma}_t^{\varepsilon})^{-1}\bar{\Gamma}_s^{\varepsilon}\right]\\
					&\leq\exp\left((1+\bar{\sigma}^2)\int_{0}^{T}u(s)ds\right).
				\end{aligned}
			\end{equation}
			Then, by equality \eqref{eq18}-\eqref{eq19}, we deduce
			\[
			\begin{aligned}
				\hat{Y}_t
				&\leq (\tilde{\Gamma}_t^{\varepsilon})^{-1}\hat{\mathbb{E}}_t^{\widetilde{G}}\left[\int_t^T\left(|m_s|+\bar{\sigma}^2|n_s|\right) \tilde{\Gamma}_s^{\varepsilon} d s\right. \left.+\int_t^T\left(|m_s^{\varepsilon}|+\bar{\sigma}^2|n_s^{\varepsilon}|\right) \tilde{\Gamma}_s^{\varepsilon} ds\right]\\
				& \leq \int_{t}^{T} C(\bar{\sigma})v(s)\phi(\frac{2L}{n-L})(\tilde{\Gamma}_t^{\varepsilon})^{-1}\hat{\mathbb{E}}_t^{\widetilde{G}}[\tilde{\Gamma}_s^{\varepsilon}] d s+\int_{t}^{T}C(\bar{\sigma}) \varepsilon (u(s)+nv(s))(\tilde{\Gamma}_t^{\varepsilon})^{-1}\hat{\mathbb{E}}_t^{\widetilde{G}}[\tilde{\Gamma}_s^{\varepsilon}] ds\\
				&\leq C(\bar{\sigma},\Lambda(u, v))\exp\left(C(\bar{\sigma},\Lambda(u, v))\right)\phi(\frac{2L}{n-L})+\varepsilon C(n,\bar{\sigma},\Lambda(u, v))\exp\left(C(\bar{\sigma},\Lambda(u, v))\right).
			\end{aligned}
			\]
			Therefore, the second term of the above inequality approaches 0 as $\varepsilon \rightarrow 0$. Hence,
			\[
			\hat{Y}_t \leq  C\left(\bar{\sigma},\Lambda(u, v)\right)\phi(\frac{2L}{n-L}),
			\]
			which ends the proof.	\end{proof}
	\end{lemma}
	
	\begin{remark}
		For any $\eta \in M_{G}^1(0,T)$, there is  
		\begin{equation*}
			\int_{0}^{T} \eta_s d\langle B\rangle_s \leq \bar{\sigma}^2 \int_{0}^{T}
			\eta_s ds.  
		\end{equation*}
	\end{remark}
	
	\begin{lemma}
		Let $\xi \in L_{G}^{\beta}(\Omega_{T})$ for some $\beta>2$. Assume that
		conditions \textnormal{(H1)}-\textnormal{(H4)} hold, and the solution to the $G $-BSDE \eqref{eq5}\eqref{eq6} is given by $\left(\bar{Y}_n,\bar{Z}_n,\bar{K}_n
		\right)$, $\left(\underline{Y}_n,\underline{Z}_n,\underline{K}_n \right)$ in $ 
		\mathfrak{S}_{G}^{2}(0, T)$ respectively. Then there
		exists a constant $C(L,\bar{\sigma},\underline{\sigma}) > 0$ such that  
		\begin{equation*}
			||\bar{Z}^n||_{H_{G}^2(0,T)}^{2} \leq C(L,\bar{\sigma},\underline{\sigma})(1+\Lambda(u,v))\left\{||\bar{Y}^n||_{S_{G}^2(0,T)}^{2}+||\bar{Y}^n||_{S_{G}^2(0,T)}\left(\left\|\int_{0}^{T}|h_0(s)|ds\right\|_{L_{G}^2(0,T)}+\Lambda(u,v)\right )\right\} ,
		\end{equation*}
		$ ||\bar{K}_T^n|_{L_{G}^2(0,T)}^{2} \leq C(L,\bar{\sigma},\underline{\sigma})(1+\Lambda(u,v))^2\left\{||\bar{Y}^n||_{S_{G}^2(0,T)}^{2}+\left\|\int_{0}^{T}|h_0(s)|ds\right\|_{L_{G}^2(0,T)}^2+\Lambda(u,v)\right\},  $
		\vspace{1em}  \newline
		where $|h_0(s)|=|f_0(s)|+|g_0(s)|$. Moreover, $\underline{Z}^n$ and $\bar{K}_T^n$ can be obtained in the same way. 
		
		\begin{proof}
			Applying the $ G \text {-Itô's} $ formula to $|\bar{Y}_t^n|^2$, we have
			\begin{equation}\label{eq10}
				\begin{aligned}
					|\bar{Y}_0^n|^2+\int_{0}^{T}|\bar{Z}_s^n|^2 d\langle B\rangle_s=&|\xi|^2+\int_{0}^{T}2\bar{Y}_s^n\left(\bar{f}_n(s,\bar{Y}_s^n,\bar{Z}_s^n)+f_0(s)\right)ds-\int_{0}^{T}2\bar{Y}_s^n\bar{Z}_s^n d B_s-\int_{0}^{T}2\bar{Y}_s^nd \bar{K}_s^n\\
					&+\int_{0}^{T}2\bar{Y}_s^n\left(\bar{g}_n(s,\bar{Y}_s^n,\bar{Z}_s^n)+g_0(s)\right)d\langle B\rangle_s.
				\end{aligned}
			\end{equation}
			From the property \textnormal{(\romannumeral 1)} in Lemma \textnormal{3.2}, it can be seen that for $ \varphi=f$, $g$, we can get 
			\[
			\left|2\bar{Y}_s^n\left(\bar{\varphi}_n(s,\bar{Y}_s^n,\bar{Z}_s^n)+\varphi_0(s)\right)\right| 
			\leq 2Lv(s)|\bar{Y}_s^n|+2Lu(s)|\bar{Y}_s^n|^2+2Lv(s)|\bar{Y}_s^n||\bar{Z}_s^n|+2L|\bar{Y}_s^n||\varphi_0(s)|.
			\]
			Given that $|h_0(s)|=|f_0(s)|+|g_0(s)|$, equation \eqref{eq10} can be scaled down to
			\[
			\begin{aligned}
				\int_{0}^{T}|\bar{Z}_s^n|^2d\langle B\rangle_s \leq &
				|\xi|^2+C(L,\bar{\sigma})\int_{0}^{T}|\bar{Y}_s^n||h_0(s)|\mathrm{d}s+C(L,\bar{\sigma})\int_{0}^{T}v(s)|\bar{Y}_s^n|ds\\
				&+C(L,\bar{\sigma})\int_{0}^{T}u(s)|\bar{Y}_s^n|^2d s+C(L,\bar{\sigma})\int_{0}^{T}v(s)|\bar{Y}_s^n||\bar{Z}_s^n|ds\\
				&+\left|\int_{0}^{T}2\bar{Y}_s^n\bar{Z}_s^ndB_s\right|+\left|\int_{0}^{T}2\bar{Y}_s^nd \bar{K}_s^n\right|.
			\end{aligned}
			\]
			Next, for any $ \varepsilon >0 $, we will show that estimates for each term in the above inequality, that is
			\[
			\hat{\mathbb{E}}\left[\int_{0}^{T}v(s)|\bar{Y}_s^n|ds\right]
			\leq \hat{\mathbb{E}}\left[\sup_{s\in[0,T]}|\bar{Y}_s^n|\int_{0}^{T}v(s)ds\right]
			\leq
			\int_{0}^{T}v(s)ds\hat{\mathbb{E}}\left[\sup_{s\in[0,T]}|\bar{Y}_s^n|^2\right]^\frac{1}{2},
			\]
			$$
			\begin{aligned}
				\hat{\mathbb{E}}\left[\int_{0}^{T}v(s)|\bar{Y}_s^n||\bar{Z}_s^n|ds\right] 
				&\leq \hat{\mathbb{E}}\left[\sup_{s\in[0,T]}|\bar{Y}_s^n|\int_{0}^{T}v(s)|\bar{Z}_s^nds\right]\\
				&\leq \hat{\mathbb{E}}\left[ \varepsilon\left(\int_{0}^{T}v(s)|\bar{Z}_s^n|ds\right )^2+\frac{1}{4\varepsilon}\sup_{s\in[0,T]}|\bar{Y}_s^n|^2\right]
				\\
				&\leq \varepsilon\int_{0}^{T}v^2(s)ds \hat{\mathbb{E}}\left[\int_{0}^{T}|\bar{Z}_s^n|^2ds\right]+\frac{1}{4\varepsilon}\hat{\mathbb{E}}\left[\sup_{s\in[0,T]}|\bar{Y}_s^n|^2\right],
			\end{aligned}
			$$
			$$
			\hat{\mathbb{E}}\left[\left|\int_{0}^{T}2\bar{Y}_s^nd \bar{K}_s^n\right|\right] \leq 2\hat{\mathbb{E}}\left[\sup_{s\in[0,T]}|\bar{Y}_s^n||\bar{K}_T^n|\right]\leq 2\hat{\mathbb{E}}[\sup_{s\in[0,T]}|\bar{Y}_s^n|^2]^\frac{1}{2}\hat{\mathbb{E}}\left[|\bar{K}_T^n|^2\right]^\frac{1}{2}.
			$$
			Applying the BDG's inequality ( see Proposition \textnormal{2.2}), it follows that
			$$
			\begin{aligned}
				\hat{\mathbb{E}}\left[\left|\int_{0}^{T}2\bar{Y}_s^n\bar{Z}_s^ndB_s\right| \right]& \leq  C(\bar{\sigma})\hat{\mathbb{E}}\left[\left|\int_{0}^{T}|\bar{Y}_s^n|^2|\bar{Z}_s^n|^2ds\right|^{\frac{1}{2}}\right]\\
				&\leq C(\bar{\sigma})\hat{\mathbb{E}}\left[\sup_{s\in[0,T]}|\bar{Y}_s^n|\left|\int_{0}^{T}|\bar{Z}_s^n|^2ds\right|^{\frac{1}{2}}\right] \\
				& \leq C(\bar{\sigma})\left\{ \varepsilon\hat{\mathbb{E}}\left[\int_{0}^{T}|\bar{Z}_s^n|^2ds\right]+\frac{1}{4\varepsilon}\hat{\mathbb{E}}\left[\sup_{s\in[0,T]}|\bar{Y}_s^n|^2\right]\right\}.
			\end{aligned}
			$$
			And note that \[
			\underline{\sigma}^2\hat{\mathbb{E}}\left[\int_{0}^{T}|\bar{Z}_s^n|^2 ds\right]
			\leq 
			\hat{\mathbb{E}}\left[\int_{0}^{T}|\bar{Z}_s^n|^2 d\langle B\rangle_s\right].
			\]
			Then, by combining the results of the above inequalities and considering the non-degeneracy of $ G $, we conclude that
			\begin{equation}\label{eq11}
				\begin{aligned}
					\underline{\sigma}^2\hat{\mathbb{E}}\left[\int_{0}^{T}|\bar{Z}_s^n|^2 ds\right]
					\leq& 
					C(\bar{\sigma},L)\left(1+\frac{1}{\varepsilon}+\int_{0}^{T}u(s)ds\right)\left\|\bar{Y}^n\right\|_{S_G^2}^2\\
					&+C(\bar{\sigma},L)\left\{\hat{\mathbb{E}}\left[\left(\int_{0}^{T}|h_0(s)|ds\right)^2\right]^\frac{1}{2}+\hat{\mathbb{E}}\left[|\bar{K}_T^n|^2\right]^\frac{1}{2}+\int_{0}^{T}v(s)ds\right\}\left\|\bar{Y}^n\right\|_{S_G^2}\\
					&+\varepsilon C(\bar{\sigma},L)\left(1+\int_{0}^{T}v^2(s)ds\right)\hat{\mathbb{E}}\left[\int_{0}^{T}|\bar{Z}_s^n|^2 ds\right].
				\end{aligned}
			\end{equation}
			On the other hand, we proceed to estimate $\bar{K}_{T}^{n}$ in the same way. And the $ G $-BSDE \eqref{eq6} is converted into
			\begin{equation}\label{eq20}
				\begin{aligned}
					\bar{K}_{T}^{n}=&\xi-\bar{Y}_{0}^{n}+\int_{0}^{T}\left[\bar{f}_{n}\left(s, \bar{Y}_{s}^{n}, \bar{Z}_{s}^{n}\right)+f_{0}(s)\right]d s\\
					&+\int_{0}^{T}\left[\bar{g}_{n}\left(s, \bar{Y}_{s}^{n}, \bar{Z}_{s}^{n}\right)+g_{0}(s)\right] d\left\langle B\right\rangle_{s}
					-\int_{0}^{T} \bar{Z}_{s}^{n} dB_{s}.
				\end{aligned}	
			\end{equation}
			By a simple calculation using the same method, we obtain
			\begin{equation}\label{eq14}
				\begin{aligned}
					\hat{\mathbb{E}}[|\bar{K}_{T}^{n}|^2]\leq & C(L,\bar{\sigma})\left(1+\left(\int_{0}^{T}u(s)ds\right )^2\right)\left\|\bar{Y}^n\right\|_{S_G^2}^2\\
					&+C(L,\bar{\sigma})\left(1+\int_{0}^{T}v^2(s)ds \right)\hat{\mathbb{E}}\left[\int_{0}^{T}|\bar{Z}_s^n|^2 ds\right]\\
					&+C(L,\bar{\sigma})\left(\int_{0}^{T}v(s)ds\right)^2+C(L,\bar{\sigma})\hat{\mathbb{E}}\left[\left(\int_{0}^{T}|h_0(s)|ds\right)^2\right]
				\end{aligned}
			\end{equation}
			and
			\begin{equation}\label{eq12}
				\begin{aligned}
					\hat{\mathbb{E}}\left[|\bar{K}_{T}^{n}|^2\right]^\frac{1}{2}\left\|\bar{Y}^n\right\|_{S_G^2}
					\leq & C(L,\bar{\sigma})\left(1+\frac{1}{4\varepsilon}+\int_{0}^{T}u(s)ds \right)\left\|\bar{Y}^n\right\|_{S_G^2}^2\\
					&+\varepsilon C(L,\bar{\sigma})\left(1+\int_{0}^{T}v^2(s)ds\right)\hat{\mathbb{E}}\left[\int_{0}^{T}|\bar{Z}_s^n|^2 ds\right]\\
					&+C(L,\bar{\sigma})\left\{\int_{0}^{T}v(s)ds+\hat{\mathbb{E}}\left[\left(\int_{0}^{T}|h_0(s)|ds\right)^2\right]^\frac{1}{2}\right\}
					\left\|\bar{Y}^n\right\|_{S_G^2}.
				\end{aligned}
			\end{equation}
			Consequently, choosing an appropriate $\varepsilon >0 $, inequalities \eqref{eq11} and \eqref{eq12} imply that
			\begin{equation}\label{eq13}
				\begin{aligned}
					\hat{\mathbb{E}}\left[\int_{0}^{T}|\bar{Z}_s^n|^2 ds\right] \leq& C(L,\bar{\sigma},\underline{\sigma})(1+\Lambda(u,v))\left\{{\hat{\mathbb{E}}\left[\sup_{s\in[0,T]}|\bar{Y}_s^n|^2\right]}\right.\\
					&\left.+{\hat{\mathbb{E}}\left[\sup_{s\in[0,T]}|\bar{Y}_s^n|^2\right]^\frac{1}{2}\left(\hat{\mathbb{E}}\left[\left(\int_{0}^{T}|h_0(s)|ds\right)^2\right]^\frac{1}{2}+\int_{0}^{T}v(s)ds\right)}
					\right\}.
				\end{aligned}
			\end{equation}
			And recalling inequality \eqref{eq14}, we deduce that
			\begin{equation}\label{eq16}
				\begin{aligned}
					\hat{\mathbb{E}}[|\bar{K}_{T}^{n}|^2] \leq& C(L,\bar{\sigma},\underline{\sigma})(1+\Lambda(u,v))^2\left\{{\hat{\mathbb{E}}\left[\sup_{s\in[0,T]}|\bar{Y}_s^n|^2\right]}\right.\\
					&\left.+{\hat{\mathbb{E}}\left[\left(\int_{0}^{T}|h_0(s)|ds\right)^2\right]+\left(\int_{0}^{T}v(s)ds\right)^2}
					\right\},
				\end{aligned}
			\end{equation}
			where $C(L,\bar{\sigma},\underline{\sigma})$ is a constant depending on $L$, $\bar{\sigma}$, and $\underline{\sigma}$. Similarly, $\underline{Z}^n$ and $\bar{K}_T^n$ can be proven to be the same, the proof is complete.
		\end{proof}
	\end{lemma}
	
	In order to ensure the convergence of $ \left\{\bar{Z}^n\right\}_{n\in \mathbb{N}} $ and $ \left\{\underline{Z}^n\right\}_{n\in \mathbb{N}} $, we need to establish some estimates for them. Note that the Lipschitz coefficient of approximating sequences ( see \eqref{eq5} \eqref{eq6}) with respect to $ z $ is $ nv(s) $, which depends on $ n $. This may lead to a blow-up in the upper bounds of $\bar{Z}^n - \bar{Z}^m$ and $\underline{Z}^n - \underline{Z}^m$, which influences on the convergence. Therefore, directly adopting the scaling idea from Lemma 3.7 is not applicable in this case. To solve this, by eliminating $ n $ during the scaling process and fully using Lemma 3.2, we provide a more precise estimate for $\bar{Z}^n - \bar{Z}^m$ and $\underline{Z}^n - \underline{Z}^m$ in $M_G^2(0,T)$, as detailed in the following Lemma.
	
	\begin{lemma}
		Set $\hat{Y}=\bar{Y}^n-\bar{Y}^m,\hat{Z}=\bar{Z}^n-\bar{Z}^m$ and $\xi \in L_{G}^{\beta}(\Omega_{T})$ for some $\beta>2$. Assume that conditions \textnormal{(H1)}-\textnormal{(H4)} hold, and the solution to the $G $-BSDE \eqref{eq6} corresponding to the generator $(\bar{f}_l+f_0(s), \bar{g}_l+g_0(s)) $ is given by $\left(\bar{Y}^l,\bar{Z}^l,\bar{K}^l\right) $ in $
		\mathfrak{S}_{G}^{2}(0, T)$, where $l=n$, $ m > L$. Then there exists a constant $C(L,\bar{\sigma},\underline{\sigma},\Lambda(u,v)) > 0$ such that 
		\begin{equation*}
			\begin{aligned}
				||\hat{Z}||_{H_{G}^2(0,T)}^{2} 
				\leq& C(L,\bar{\sigma},\underline{\sigma})(1+\Lambda(u,v))^2 
				\Bigl\{ 
				||\hat{Y}||_{S_{G}^2(0,T)}^{2} 
				+ ||\hat{Y}||_{S_{G}^2(0,T)} \Bigl( 
				||\bar{Y}^n||_{S_{G}^2(0,T)} + ||\bar{Y}^m||_{S_{G}^2(0,T)} 
				\Bigr. \\
				&\Bigl. 
				+ \left\|\int_{0}^{T} |h_0(s)| ds \right\|_{L_{G}^2(0,T)} 
				+ \Lambda(u,v) 
				\Bigr) 
				\Bigr\},
			\end{aligned}
		\end{equation*}
		where $C(L,\bar{\sigma},\underline{\sigma})$ is a constant depending on $L$, $\bar{\sigma}$, $\underline{\sigma}$. And the estimate for $\underline{Z}^n_t-\underline{Z}^m_t$ can be proved in the similar way. 
		
		\begin{proof}
			The proof is generally similar to that of Lemma \textnormal{3.7}, we will focus on the different parts. Applying the $ G $$  \text {-Itô's formula } $to $|\hat{Y}_t|^2=|\bar{Y}^n_t-\bar{Y}^m_t|^2$ yields that
			\begin{equation}\label{15}
				\begin{aligned}
					|\hat{Y}_0|^2+\int_{0}^{T}|\hat{Z}_s|^2 d\langle B\rangle_s =|\xi|^2+\int_{0}^{T}2\hat{Y}_s\hat{f}(s)ds+\int_{0}^{T}2\hat{Y}_s\hat{g}(s)d\langle B\rangle_s-\int_{0}^{T}2\hat{Y}_s\hat{Z}_sd B_s-\int_{0}^{T}2\hat{Y}_sd\hat{K}_s,
				\end{aligned}
			\end{equation}
			where
			\[
			\hat{f}(s)=\bar{f}_n(s,\bar{Y}_s^n,\bar{Z}_s^n)-\bar{f}_m(s,\bar{Y}_s^m,\bar{Z}_s^m),\text{ }
			\hat{g}(s)=\bar{g}_n(s,\bar{Y}_s^n,\bar{Z}_s^n)-\bar{g}_m(s,\bar{Y}_s^m,\bar{Z}_s^m).
			\]
			Set $\hat{h}_s=|\bar{f}_n(s,\bar{Y}_s^m,\bar{Z}_s^n)-\bar{f}_m(s,\bar{Y}_s^m,\bar{Z}_s^m)|+|\bar{g}_n(s,\bar{Y}_s^m,\bar{Z}_s^n)-\bar{g}_m(s,\bar{Y}_s^m,\bar{Z}_s^m)|$. Recalling (H2), we have
			\begin{equation}\label{eq29}
			\begin{aligned}
				|\hat{\varphi}(s)|&=|\bar{\varphi}_n(s,\bar{Y}_s^n,\bar{Z}_s^n)-\bar{\varphi}_n(s,\bar{Y}_s^m,\bar{Z}_s^n)+\bar{\varphi}_n(s,\bar{Y}_s^m,\bar{Z}_s^n)-\bar{\varphi}_m(s,\bar{Y}_s^m,\bar{Z}_s^m)|\\
				& \leq u(s)|\bar{Y}_s^n-\bar{Y}_s^m|+|\bar{\varphi}_n(s,\bar{Y}_s^m,\bar{Z}_s^n)-\bar{\varphi}_m(s,\bar{Y}_s^m,\bar{Z}_s^m)|\\
				& \leq u(s)|\hat{Y}_s|+\hat{h}_s,
			\end{aligned}
		  \end{equation}
			where $ \varphi=f$, $g $. Then, by Lemma \textnormal{3.2}, it follows that
			$$
			\begin{aligned}
				2|\hat{Y}_s\hat{\varphi}(s)|
				\leq 2|\hat{Y}_s|(u(s)|\hat{Y}_s|+\hat{h}_s)=2u(s)|\hat{Y}_s|^2+2\hat{h}_s|\hat{Y}_s|.
			\end{aligned}
			$$
			Based on the BDG's inequality, it holds that
			$$
			\begin{aligned}
				\hat{\mathbb{E}}\left[\left|\int_{0}^{T}2\hat{Y}_s\hat{Z}_sdB_s\right| \right]& \leq  C(\bar{\sigma})\hat{\mathbb{E}}\left[\left|\int_{0}^{T}|\hat{Y}_s|^2|\hat{Z}_s|^2ds\right|^{\frac{1}{2}}\right]\\
				&\leq C(\bar{\sigma})\hat{\mathbb{E}}\left[|\sup_{s\in[0,T]}|\hat{Y}_s|\left|\int_{0}^{T}|\hat{Z}_s|^2ds\right|^{\frac{1}{2}}\right] \\
				& \leq C(\bar{\sigma})\left\{ \varepsilon^{\prime}\hat{\mathbb{E}}\left[\int_{0}^{T}|\hat{Z}_s|^2ds\right]+\frac{1}{4\varepsilon^{\prime}}\hat{\mathbb{E}}\left[\sup_{s\in[0,T]}|\hat{Y}_s|^2\right]\right\}.
			\end{aligned}
			$$
			Then, by simple calculation, we obtain that                                                                                                                                                         
			\begin{equation}\label{eq28}
				\begin{aligned}
					\underline{\sigma}^2\hat{\mathbb{E}}\left[\int_{0}^{T}|\hat{Z}_s|^2 ds\right]
					\leq& C(\bar{\sigma},L)\left(\frac{1}{\varepsilon^{\prime}}+\int_{0}^{T}u(s)ds\right)||\hat{Y}||_{S_G^2}^2+\varepsilon^{\prime}C(\bar{\sigma},L)\hat{\mathbb{E}}\left[\int_{0}^{T}|\hat{Z}_s|^2 ds\right]\\
					&+C(\bar{\sigma},L)||\hat{Y}||_{S_G^2}\left(\hat{\mathbb{E}}\left[\left(\int_{0}^{T}|\hat{h}_s|ds\right)^2\right]^\frac{1}{2}+\hat{\mathbb{E}}\left[(|\bar{K}_T^n|+|\bar{K}_T^m|)^2\right]^\frac{1}{2}\right).
				\end{aligned}
			\end{equation}
			Following this, we proceed to scale the second term of inequality \eqref{eq28}. Since the sequence of generator $ \bar{\varphi}_n $ can be uniformly controlled, for each $ n$, $m>L $, we have
			$$
			\hat{h}_s \leq 4L(u(s)|\bar{Y}_s^m|+v(s)|\bar{Z}_s^m|+v(s)|\bar{Z}_s^n|+v(s)).
			$$
			Set $ C(L)=C(L,\bar{\sigma},\underline{\sigma}).$ Recalling the estimation results from the previous inequalities \eqref{eq13} and \eqref{eq16} in Lemma \textnormal{3.7}, we derive  
			$$
			\begin{aligned}
				\hat{\mathbb{E}}\left[\left(\int_{0}^{T}|\hat{h}_s|ds\right)^2\right]^\frac{1}{2} 
				\leq& 4L\hat{\mathbb{E}}\left[\left(\int_{0}^{T}(u(s)|\bar{Y}_s^m|+v(s)|\bar{Z}_s^m|+v(s)|\bar{Z}_s^n|+v(s))ds\right)^2\right]\\
				\leq& C(L)\left\{\int_{0}^{T}u(s)ds\left\|\bar{Y}^m\right\|_{S_G^2}+\left(\int_{0}^{T}v^2(s)ds\right)^\frac{1}{2}(\left\|\bar{Z}^m\right\|_{H_G^2}+\left\|\bar{Z}^n\right\|_{H_G^2})+\int_{0}^{T}v(s)ds\right\}\\
				\leq& C(L)(1+\Lambda(u,v))^2\left\{\left\|\bar{Y}^m\right\|_{S_G^2}+\left\|\bar{Y}^n\right\|_{S_G^2}+\hat{\mathbb{E}}\left[\left(\int_{0}^{T}|h_0(s)|ds\right)^2\right]^\frac{1}{2}+\int_{0}^{T}v(s)ds\right\}.
			\end{aligned}
			$$
			Thus, we can get that
			$$ 
			\begin{aligned}
				\hat{\mathbb{E}}\left[\left(\int_{0}^{T}|\hat{h}_s|ds\right)^2\right]^\frac{1}{2} + \hat{\mathbb{E}}\left[(|\bar{K}_T^n|+|\bar{K}_T^m|)^2\right]^\frac{1}{2} \leq& C(L)(1+\Lambda(u,v))^2 \left\{{\hat{\mathbb{E}}\left[\sup_{s\in[0,T]}|\bar{Y}_s^n|^2\right]^\frac{1}{2}+\hat{\mathbb{E}}\left[\sup_{s\in[0,T]}|\bar{Y}_s^m|^2\right]^\frac{1}{2}} \right.\\
				&\left.+ {\hat{\mathbb{E}}\left[\left(\int_{0}^{T}|h_0(s)|ds\right)^2\right]^\frac{1}{2} + \int_{0}^{T}v(s)ds} \right\}.
			\end{aligned}
			$$
			In conclusion, by selecting an appropriate $\varepsilon^{\prime}>0$, from \eqref{eq28}, it follows from \eqref{eq28} that
			$$
			\begin{aligned}
				\hat{\mathbb{E}}\left[\int_{0}^{T}|\hat{Z}_s|^2 ds\right]
				\leq& C(L)(1+\Lambda(u,v))\hat{\mathbb{E}}\left[\sup_{s\in[0,T]}|\hat{Y}_s|^2\right]\\
				&+C(L)(1+\Lambda(u,v))^2\hat{\mathbb{E}}\left[\sup_{s\in[0,T]}|\hat{Y}_s|^2\right]^{\frac{1}{2}}\left\{{\hat{\mathbb{E}}\left[\sup_{s\in[0,T]}|\bar{Y}_s^n|^2\right]^\frac{1}{2}+\hat{\mathbb{E}}\left[\sup_{s\in[0,T]}|\bar{Y}_s^m|^2\right]^\frac{1}{2}}\right.\\
				&\left.{+\hat{\mathbb{E}}\left[\left(\int_{0}^{T}|h_0(s)|ds\right)^2\right]^\frac{1}{2}+\int_{0}^{T}v(s)ds}\right\}.
			\end{aligned}
			$$
			With this, the proof is complete.
		\end{proof}
	\end{lemma}
	
	The Lemmas have been verified. Next we will apply the conclusions of the Lemmas to prove the existence and uniqueness of solutions for $G $-BSDE \eqref{eq1}. And the main results of this
	paper are as follows.
	\begin{theorem}
		Assume that the conditions \textnormal{(H1)}-\textnormal{(H4)} hold and $\xi \in L_{G}^{\beta}(\Omega_T)$ for some $\beta >2$. Then the $G $-BSDE \eqref{eq1} admits a
		unique solution $(Y,Z,K)\in \mathfrak{S}_{G}^{2}(0, T)$.  
		
		\begin{proof} The proof can be divided into two parts.
			\\
			\textbf{Part I: The uniqueness of the solution}
			
			Suppose that $(Y^1,Z^1,K^1)$ and $(Y^2,Z^2,K^2)\in \mathfrak{S}_{G}^{2}(0, T)$ are both solutions of the $ G $-BSDE \eqref{eq1}. By Lemma \textnormal{3.4} and Theorem 2.5
			we deduce that for each $ n>L$,
			$$
			\underline{Y}_t^n \leq Y_t^i \leq \bar{Y}_t^n,\quad i=1,2.
			$$
			Following this, based on the estimation of $\underline{Y}_t^n-\bar{Y}_t^n$ in Lemma \textnormal{3.5}, we have
			\[
			|Y_t^1-Y_t^2| \leq |\bar{Y}_t^n-%
			\underline{Y}_t^n| \leq C\left(\bar{\sigma}, \underline{\sigma},\Lambda(u,
			v)\right)\phi(\frac{2L}{n-L}),\quad \forall n>L,\text{ } \forall t \in
			[0,T].
			\]
			Since the non-decreasing property of the function $\phi(\cdot)$ and $Y_t^i$ is continuous process, we have $Y_t^1=Y_t^2$ as $ n \longrightarrow \infty$.
			
			On the other hand, by applying $ G \text {-Itô's formula} $ to $|Y_t^1-Y_t^2|^2$, the
			uniqueness of the solution $(Z,K)$ also can be proven. Set $\hat{Y}_t=Y_t^1-Y_t^2$, we can get that
			$$
			|\hat{Y}_0|^2+\int_{0}^{T}|\hat{Z}_s|^2 d\langle B\rangle_s=|\hat{Y}
			_T|^2+\int_{0}^{T}2\hat{Y}_s\hat{f}(s)ds +\int_{0}^{T}2\hat{Y}_s\hat{g}
			(s)d\langle B\rangle_s-\int_{0}^{T}2\hat{Y}_s\hat{Z}_sd B_s-\int_{0}^{T}2
			\hat{Y}_sd \hat{K}_s,
			$$
			where $\hat{f}(s)=f(s,Y_s^1,Z_s^1)-f(s,Y_s^2,Z_s^2)$ and $\hat{g}
			(s)=g(s,Y_s^1,Z_s^1)-g(s,Y_s^2,Z_s^2)$. Note that $\hat{Y}_t=0 $ for any $ t
			\in [0,T]$, so we can conclude that	$Z_t^1=Z_t^2$ and $K_t^1=K_t^2$, which shows the uniqueness of the solution to $ G $-BSDE \eqref{eq1}.
			\\
			\textbf{Part II: The existence of the
				solution}
			
			The main idea of the proof is to make full use of the convergence of the solution $ (\bar{Y}^n, \bar{Z}^n, \bar{K}^n) $ for the $ G $-BSDE \eqref{eq6}. The entire process will be carried out in three steps.
			\\
			Step \textnormal{1}: the uniform boundedness of $(\bar{Y}^n,\bar{Z}^n,\bar{K}^n)$
			\\
			Based on Lemma \textnormal{3.4} and Lemma \textnormal{3.7}, there exists a constant $C\left(L,\bar{\sigma}, 
			\underline{\sigma},\Lambda(u, v)\right) >0$ (simply noted as C) such that for each $ n $, we derive that
			$$
			||\bar{Y}^n||_{S_{G}^2(0,T)}\leq C,\text{ }||\bar{Z}^n||_{M_{G}^2(0,T)}\leq C,\text{ }||
			\bar{K}_T^n||_{L_{G}^2(0,T)}\leq C.
			$$
			Step \textnormal{2}: the convergence of $(\bar{Y}^n,\bar{Z}^n,\bar{K}^n)$
			
			We first show that the convergence of the sequence $\left\{\bar{Y}^n\right\}_{n\in N}$. From Lemma \textnormal{3.4} and a priori estimate for $|\bar{Y}^n - \underline{Y}^n|$ in Lemma \textnormal{3.5}, for any 
			$n, m > L$, we have
			$$
			\begin{aligned}
				||\bar{Y}^n-\bar{Y}^m||_{S_{G}^2(0,T)}^2 &\leq ||\bar{Y%
				}^{n\wedge m}-\underline{Y}^{n\wedge m}||_{S_{G}^2(0,T)}^2=\hat{\mathbb{E}}\left[%
				\sup_{s\in[0,T]}|\bar{Y}^{n\wedge m}_s-\underline{Y}^{n\wedge m}_s|^2\right]%
				\\
				& \leq C\left(\bar{\sigma}, \underline{\sigma},\Lambda(u,
				v)\right)\phi^2(\frac{2L}{n\wedge m-L}).
			\end{aligned}
			$$
			Sending $n$, $m\rightarrow \infty $, it holds that
			\[
			\lim_{n \to \infty}||\bar{Y}^n-\bar{Y}^m||_{S_{G}^2(0,T)}^2 =0.
			\] 
			Then it is clear that $\left\{\bar{Y}^n\right\}_{n \in N}$ is a Cauchy sequence, i.e. there exists a process $Y \in S_{G}^2(0,T)$ such that
			\begin{equation}\label{eq31}
				\lim_{n\to \infty}||\bar{Y}^n-Y||_{S_{G}^2(0,T)}^2 =0.
			\end{equation}
			
			Next, we show that the convergence of $\{%
			\bar{Z}^n\}_{n \in N}$. Based on a priori estimate of $\hat{Z}$ in Lemma \textnormal{3.8}, we can derive that
			$$
			\begin{aligned}
				\lim_{n,m \to \infty}||\bar{Z}^n-\bar{Z}^m||_{H_{G}^2(0,T)}^2 = 0.
			\end{aligned}
			$$
			Then $\left\{\bar{Z}^n\right\}_{n \in N}$ is also a Cauchy
			sequence which indicates that  there also exists a process 
			$Z \in H_{G}^2(0,T)$ such that
			\begin{equation}\label{eq32}
				\lim_{n,m \to \infty}||\bar{Z}^n-Z||_{H_{G}^2(0,T)}^2=0.
			\end{equation}
			
			Finally, we show that the convergence of $\bar{K}_T^n$ in $L_G^2(\Omega_T)$. By simplifying the $ G $-BSDE \eqref{eq1} and \eqref{eq6},
			we can obtain that
			\[
			K_t=Y_t-Y_0+\int_{0}^{T}f(s,Y_s,Z_s)ds+%
			\int_{0}^{T}g(s,Y_s,Z_s)d\langle B \rangle
			_s-\int_{0}^{T}Z_sdB_s
			\]
			\[
			\bar{K}_t^n=\bar{Y}_t^n-\bar{Y}_0^n+\int_{0}^{T}(\bar{f}_n(s,\bar{Y}_s^n,\bar{Z}_s^n)+f_0(s))ds+%
			\int_{0}^{T}(\bar{g}_n(s,\bar{Y}_s^n,\bar{Z}_s^n)+g_0(s))d\langle B \rangle
			_s-\int_{0}^{T}\bar{Z}_s^ndB_s
			\]
			\\
			Combining this with the convergence of $\left\{\bar{Y}^n\right\}_{n \in N}$ and $\{\bar{Z}^n\}_{n\in N}$, we only need to prove 
			\begin{equation}\label{17}
				\lim_{n \to\infty} \hat{\mathbb{E}}\left[(\int_{0}^{T}\bar{f}_n(s,\bar{Y}_s^n,\bar{Z}_s^n)+f_0(s)-f(s,Y_s,Z_s)ds)^2\right]=0.
			\end{equation}
			\\
			Recalling the assumption \textnormal{(H2)} and Lemma \textnormal{3.2 (\romannumeral 5)}, since
			$$
			\begin{aligned}
				&|\bar{f}_n(s,\bar{Y}_s^n,\bar{Z}
				_s^n)+f_0(s)-f(s,Y_s,Z_s)| \\
				\leq& |\bar{f}_n(s,\bar{Y}_s^n,\bar{Z}
				_s^n)+f_0(s)-f(s,\bar{Y}_s^n,\bar{Z}_s^n)| +|f(s,\bar{Y}_s^n,\bar{Z}%
				_s^n)-f(s,Y_s,\bar{Z}_s^n)|+|f(s,\bar{Y}_s^n,\bar{Z}_s^n)-f(s,Y_s,Z_s)|\\
				\leq& v(s)\phi(\frac{2L}{n-L})+u(s)|\bar{Y}_s^n-Y_s|+v(s)\phi(|\bar{Z}
				_s^n-Z_s|),
			\end{aligned}
			$$
			we have
			$$
			\begin{aligned}
				&\hat{\mathbb{E}}\left[\left(\int_{0}^{T}\bar{f}_n(s,\bar{Y}_s^n,\bar{Z}_s^n)+f_0(s)-f(s,Y_s,Z_s)ds\right)^2\right]\\
				\leq  &3\left\{\hat{\mathbb{E}}\left[\left(\int_{0}^{T}v(s)\phi(\frac{2L}{n-L})ds\right )^2\right]+\hat{\mathbb{E}}\left[\left(\int_{0}^{T}u(s)|\bar{Y}_s^n-Y_s|ds\right )^2\right]+\hat{\mathbb{E}}\left[\left(\int_{0}^{T}v(s)\phi(|\bar{Z}_s^n-Z_s|)ds\right )^2\right]\right\}\\
				:=&3(I_1+I_2+I_3).
			\end{aligned}
			$$
			Then, we will proceed to estimate $I_1, I_2, I_3$ respectively. It is easy to verify that
			$$
			I_1\leq \Lambda(u,v)\phi^2(\frac{2L}{n-L}),
			$$
			$$
			I_2\leq \hat{\mathbb{E}}\left[\sup_{s\in[0,T]}|\bar{Y}_s^n-Y_s|^2\left(\int_{0}^{T}u(s)ds\right)^2\right]\leq\Lambda^2(u,v)\hat{\mathbb{E}}\left[\sup_{s\in[0,T]}|\bar{Y}_s^n-Y_s|^2\right].
			$$
			\\
			Based on the non-decreasing property of $\phi$ and the convergence of $\left\{\bar{Y}^n\right\}_{n \in N}$ (see \eqref{eq31}), we deduce that $I_1 + I_2 \to 0$ as $n\rightarrow \infty$.
			And since $\phi$ is a uniformly continuous function, it can be obtained that for any $\varepsilon > 0$, there exists a constant $\delta > 0$
			such that when $|\bar{Z}_s^n - Z_s| < \delta$, we have $\phi(|\bar{Z}_s^n -
			Z_s|) < \varepsilon$. Therefore, by taking the constant 
			$\delta$ as given above, it holds that
			\begin{equation}\label{eq30}
			\begin{aligned}
				I_3&\leq 2\hat{\mathbb{E}}\left[\left(\int_{0}^{T}v(s)\phi(|\bar{Z}_s^n-Z_s|)\mathbf{I}_{\left\{|\bar{Z}_s^n-Z_s|\leq \delta\right\}}ds\right )^2\right]+2\hat{E
				}\left[\left(\int_{0}^{T}v(s)\phi(|\bar{Z}_s^n-Z_s|)\mathbf{I}_{\left\{|\bar{Z}_s^n-Z_s|>\delta\right\}}ds\right )^2\right]\\
				&\leq 2\varepsilon^2\Lambda(u,v)+2\hat{\mathbb{E}}\left[\int_{0}^{T}v^2(s)ds\int_{0}^{T}\phi^2(|\bar{Z}_s^n-Z_s|)\mathbf{I}_{\left\{|\bar{Z}_s^n-Z_s|>\delta\right\}}ds\right]\\
				&\leq2\varepsilon^2\Lambda(u,v)+C\Lambda(u,v)\hat{\mathbb{E}}\left[\int_{0}^{T}(1+|\bar{Z}_s^n|^2+|Z_s|^2)\mathbf{I}_{\left\{|\bar{Z}_s^n-Z_s|>\delta\right\}}ds\right].
			\end{aligned}
		   \end{equation}
			Following this, applying the Marcov inequality and the convergence of $\bar{Z}^n$ in $M_G^2(0,T)$ (see \eqref{eq32}), the following equation \eqref{eq23} can be obtained. Sending $ n\rightarrow \infty $, we have
			\begin{equation}\label{eq23}
				\hat{\mathbb{E}}\left[\int_{0}^{T}\mathbf{I}_{\left\{|\bar{Z}_s^n-Z_s|> \delta\right\}}ds\right]\leq \frac{\hat{\mathbb{E}}\left[\int_{0}^{T}|\bar{Z}_s^n-Z_s|^2ds\right]}{\delta^2}\rightarrow 0.
			\end{equation}
			According to proposition \textnormal{3.8} in \cite{li2011stopping}, we can get that
			\begin{equation}\label{eq24}
				\forall \zeta  \in H_G^\alpha(0,T,R^d),\text{ }\lim_{N \to \infty}\hat{\mathbb{E}}\left[\left(\int_{0}^{T}|\zeta_s|^2\mathbf{I}_{\left\{|\zeta_s|>N\right\}}\right)^\frac{\alpha}{2}ds\right]=0.
			\end{equation}
			Thus, select a fixed $ N \in \mathbb{N} $, it follows that 
			\begin{equation}\label{eq21}
				\begin{aligned}
					\hat{\mathbb{E}}\left[\int_{0}^{T}|\bar{Z}_s^n|^2\mathbf{I}_{\left\{|\bar{Z}_s^n-Z_s|>\delta\right\}}ds\right]
					\leq &C\hat{\mathbb{E}}\left[%
					\int_{0}^{T}|\bar{Z}_s^n|^2\mathbf{I}_{\left\{|\bar{Z}_s^n-Z_s|>\delta\right\}}\mathbf{I}_{
						\left\{|\bar{Z}_s^n|>N\right\}}ds\right]\\
					&+C\hat{\mathbb{E}}\left[\int_{0}^{T}|%
					\bar{Z}_s^n|^2\mathbf{I}_{\left\{|\bar{Z}_s^n-Z_s|>\delta\right\}}\mathbf{I}_{\left\{|\bar{Z}%
						_s^n|\leq N\right\}}ds\right]\\
					\leq& C\hat{\mathbb{E}}\left[\int_{0}^{T}|\bar{Z}_s^n|^2\mathbf{I}_{\left\{|\bar{Z}_s^n|>N\right\}}ds\right]+CN^2\hat{\mathbb{E}}\left[\int_{0}^{T}\mathbf{I}_{\left\{|\bar{Z}_s^n-Z_s|>\delta\right\}}ds\right]
				\end{aligned}
			\end{equation}
			and
			\begin{equation}\label{eq22}
				\begin{aligned}
					\hat{\mathbb{E}}\left[\int_{0}^{T}|Z_s|^2\mathbf{I}_{\left\{|\bar{Z}_s^n-Z_s|>\delta\right\}}ds\right]
					\leq &C\hat{\mathbb{E}}\left[	\int_{0}^{T}|Z_s|^2\mathbf{I}_{\left\{|\bar{Z}_s^n-Z_s|>\delta\right\}}\mathbf{I}_{			\left\{|Z_s|>N\right\}}ds\right]\\
					&+C\hat{\mathbb{E}}\left[\int_{0}^{T}|		Z_s|^2\mathbf{I}_{\left\{|\bar{Z}_s^n-Z_s|>\delta\right\}}\mathbf{I}_{\left\{|Z_s|\leq N\right\}}ds\right]\\
					\leq& C\hat{\mathbb{E}}\left[\int_{0}^{T}|Z_s|^2\mathbf{I}_{\left\{|\bar{Z}_s^n|>N\right\}}ds\right]+CN^2\hat{\mathbb{E}}\left[
					\int_{0}^{T}\mathbf{I}_{\left\{|\bar{Z}_s^n-Z_s|>\delta\right\}}ds\right].
				\end{aligned}
			\end{equation}
			From the inequalities \eqref{eq30}-\eqref{eq22}, we will prove the convergence of $ I_3 $. Firstly, by fixing $N$ and sending $n \to \infty$, we deduce that
			\[
			I_3\leq C\Lambda(u,v)\left\{\varepsilon^2+\hat{\mathbb{E}}\left[\int_{0}^{T}|\bar{Z}_s^n|^2\mathbf{I}_{\left\{|\bar{Z}_s^n|>N\right\}}ds\right]+\hat{\mathbb{E}}\left[\int_{0}^{T}|Z_s|^2\mathbf{I}_{\left\{|Z_s|>N\right\}}ds\right]\right\}.
			\]
			Then, sending $ N \rightarrow \infty $ and $ \varepsilon \rightarrow 0 $, we can get $ I_3\rightarrow 0 $ by \eqref{eq24}.\\			
			In conclusion, the
			proof of equation \eqref{17} is complete. Then there exists a process $K\in
			L_G^2(\Omega_T)$ such that
			\[\lim_{n \to \infty}||\bar{K}_T^n-K_T||_{L_G^2(\Omega_T)}= 0.
			\]
			Step 3: $K_t$ is a G-martingale
			
			For all $0\leq s \leq t \leq T$, we have
			$$
			\begin{aligned}
				\hat{\mathbb{E}}\left[\left|\hat{\mathbb{E}}_s[K_t]-K_s\right|\right]
				&\leq \hat{\mathbb{E}}\left[\left|\hat{\mathbb{E}}_s[K_t]-\hat{\mathbb{E}}_s[\bar{K}_t^n]+\hat{\mathbb{E}}_s[\bar{K}_t^n]-K_s\right|\right]\\
				& \leq \hat{\mathbb{E}}\left[\hat{\mathbb{E}}_s[|K_t-%
				\bar{K}_t^n|]+|\bar{K}_s^n-K_s|  \right]\\
				&\leq \hat{\mathbb{E}}[|K_t-\bar{K}_t^n|]+\hat{\mathbb{E}}[|\bar{K}_s^n-K_s|].
			\end{aligned}
			$$
			Combining this with the convergence of $\bar{K}_t^n$ in $L_G^2(\Omega_T)$, sending $n \rightarrow \infty$, it holds that
			\[
			\hat{\mathbb{E}}[|K_t-\bar{K}_t^n|]+\hat{\mathbb{E}}[|\bar{K}_s^n-K_s|]\rightarrow 0,
			\]
			which proves $K_t$ is a $ G $-martingale. Consequently, the $ G $-BSDE \eqref{eq1} admits a unique solution in $\mathfrak{S}_{G}^{2}(0, T)$ and the main result of this paper have been proved.
		\end{proof}	
	\end{theorem}

	The Comparison Theorem presented in this paper can be demonstrated by Proposition \textnormal{2.5}.
	\begin{theorem}
		(Comparison Theorem) Suppose that $(Y^i, Z^i, K^i) $, $i=1,2 $ is the $\mathfrak{S}_{G}^{2}(0, T)$ solution to the following $ G $-BSDE corresponding to the data $(\xi^i, f^i, g^i) $ on finite interval $ [0,T] $:
		\begin{equation}
			Y_{t}^i=\xi^i +\int_{t}^{T} f^i(s,Y_{s},Z_{s})ds+\int_{t}^{T}
			g^i(s,Y_{s},Z_{s})d\langle B \rangle_{s}
			-\int_{t}^{T}Z_{s}^idB_{s}-(K_{T}^i-K_{t}^i),
		\end{equation}
		where $\xi^{i}$, $f^i$ and $g^i$ satisfy the assumptions \textnormal{(H1)-(H4)}. For any $(t,\omega,y,z) $, if $\xi^1 \leq \xi^2 $, $f^1(t,y,z)\leq f^2(t,y,z)$ and $g^1(t,y,z)\leq g^2(t,y,z) $ , then $Y_t^1\leq Y_t^2$. 
		
		\begin{proof}
			Let $(\underline{Y}^{1,n}, \underline{Z}^{1,n}, \underline{K}^{1,n}) $ be the solution of $ G $-BSDE $ \eqref{eq6} $ with respect to the data $ (\xi^1, f^1, g^1) $ for any $ n \in \mathbb{N} $. Based on Lemma \textnormal{3.2}, we deduce that
			\[
			f^2(t,y,z)\geq f^1(t,y,z) \geq \underline{f}^1_n(t,y,z)+f_0^1(t),
			\]
			\[
			g^2(t,y,z)\geq g^1(t,y,z) \geq \underline{g}^1_n(t,y,z)+g_0^1(t).
			\]
			Then recalling Theorem \textnormal{2.5}, we have $ Y_t^2\geq \underline{Y}_t^{1,n} $ for any $ t \in [0,T] $. Furthermore, according to the convergence of $Y_t^{1,n}$ in $S_{G}^{2}(0, T)$ (see \eqref{eq31}), we conclude that $ Y_t^1\leq Y_t^2$ for any $ t \in [0,T]$.
		\end{proof}
	\end{theorem}


\begin{thebibliography}{99} 
		\bibitem {buckdahn2010probabilistic}R. Buckdahn, Y. Hu, Probabilistic interpretation of a coupled system of Hamilton–Jacobi–Bellman equations, J.Evol.Equ Evolution Equations, 10(2010), 529–549.
		
		\bibitem{chen2000infinite}Z. Chen, B. Wang, Infinite time interval BSDEs and the convergence of g-martingales, J. Aust. Math. Soc., 69(2000), 187–211.
		
		\bibitem{coquet2002filtration}F. Coquet, Y. Hu, J. Mémin, S. Peng, Filtration-consistent nonlinear expectations and related g-expectations, Probab. Theory Related.Fields, 123(2002), 1–27.
		
		\bibitem{denis2011function}L. Denis, M. Hu, S. Peng, Function Spaces and Capacity Related to a Sublinear Expectation: Application to $ G $-Brownian Motion Paths, Potential Anal., 34(2011), 139–161.
		
		\bibitem{el1997backward}N. El Karoui, S. Peng, M. C. Quenez, Backward stochastic differential equations in finance, Math.Finance, 7(1997), 1–71.
		
		\bibitem{hu2014backward}M. Hu, S. Ji, S. Peng, Y. Song, Backward stochastic differential equations driven by G-Brownian motion, Stochastic Process. Appl., 124(2014), 759–784.
		
		\bibitem{hu2014comparison}M. Hu, S. Ji, S. Peng, Y. Song, Comparison Theorem, Feynman–Kac formula and Girsanov transformation for BSDEs driven by G-Brownian motion, Stochastic Process. Appl., 124(2014), 1170–1195.
		
		\bibitem{hu2020bsdes}M. Hu, B. Qu, F. Wang, BSDEs driven by G-Brownian motion with time-varying Lipschitz condition, J. Math. Anal. Appl., 491(2020), 124342.
		
		\bibitem{hu2016quasi}M. Hu, F. Wang, G. Zheng, Quasi-continuous random variables and processes under the G-expectation framework, Stochastic Process. Appl., 126(2016), 2367–2387.
		
		10\bibitem{hu2020quadraticbsde}Y. Hu, Y. Lin, A. S. Hima, Quadratic backward stochastic differential equations driven by G-Brownian motion: Discrete solutions and approximation, Stochastic Process. Appl., 128(2018), 3724–3750.
	
		\bibitem{huQuadraticGBSDEs2022}Y. Hu, S. Tang, F. Wang, Quadratic G-BSDEs with convex generators and unbounded terminal conditions, Stochastic Process. Appl., 153(2022), 363–390.
		
		\bibitem{li2011stopping}X. Li, S. Peng, Stopping times and related Itô’s calculus with G-Brownian motion, Stochastic Process. Appl., 121(2011), 1492–1508.
		
		\bibitem{lin2bsde2020}Y. Lin, Z. Ren, N. Touzi, J. Yang, Second order backward SDE with random terminal time, Electron. J. Probab., 25(2020), 1-43.
		
		\bibitem{liu2020multi}G. Liu, Multi-dimensional BSDEs driven by G-Brownian motion and related system of fully nonlinear PDEs, Stochastics, 92(2020), 659–683.
		
		\bibitem{mao1995adapted}X. Mao, Adapted solutions of backward stochastic differential equations with non-lipschitz coefficients, Stochastic Process. Appl., 58(1995), 281–292.
		
		\bibitem{pardoux1990adapted}E. Pardoux, S. Peng, Adapted solution of a backward stochastic differential equation, Syst. Control Lett, 14(1990), 55–61.
		
		\bibitem{peng2005nonlinear}S. Peng, Nonlinear expectations and nonlinear Markov chains, Chin. Ann. Math, 26(2005), 159–184.
		
		\bibitem{peng2007g}S. Peng, G-expectation, G-Brownian motion and related stochastic calculus of Itô type, Stochastic analysis and applications, Abel Symp., Springer, Berlin, 2(2007), 541-567.
		
		\bibitem{peng2008multi}S. Peng, Multi-dimensional G-Brownian Motion and Related Stochastic Calculus under G-Expectation, Stochastic Process. Appl., 118(2008), 2223–2253.
		
		\bibitem{peng2010nonlinear}S. Peng, Nonlinear Expectations and Stochastic Calculus Under Uncertainty, Springer, Berlin(2019).
	
		\bibitem{soner2012wellposedness}H. M. Soner, N. Touzi, J. Zhang, Wellposedness of second order backward SDEs, Probab. Theory Related Fields, 153(2012), 149–190.
		
		\bibitem{soner2013dual}H. M. Soner, N. Touzi, J. Zhang, Dual formulation of second order target problems, Ann. Appl. Probab., 23(2013).
		
		\bibitem{wang2021backward}F. Wang,G. Zheng, Backward stochastic differential equations driven by G-Bsrownian motion with	uniformly continuous generators, J. Theoret. Probab., 34(2021), 660–681.
		
	\end{thebibliography}
\end{document}